\documentclass[12pt]{amsart}
\usepackage{mathrsfs,amssymb,amsbsy,amsmath,booktabs,tabularx}
\usepackage[colorlinks,linkcolor=blue,citecolor=blue,pagebackref,hyperfootnotes=false]{hyperref}
\hypersetup{pdfauthor={Yude Liu, Xinbao Lu, Qiang Sun}}
\newtheorem{theorem}{Theorem}[section]
\newtheorem{lemma}[theorem]{Lemma}
\newtheorem{corollary}[theorem]{Corollary}

\makeatletter
\newcommand{\applicationsubsection}{\@startsection{subsection}{2}{\z@}%
  {.5\linespacing\@plus.7\linespacing}{.35\linespacing}{\normalfont\bfseries}}
\makeatother
\allowdisplaybreaks
\numberwithin{equation}{section}
\newcommand{\R}{\mathbb R}
\newcommand{\dd}{\,d}

\begin{document}
\begin{center}
{\Large\bfseries Brunn--Minkowski inequalities for chord integrals\\
and fractional perimeters}
\end{center}
\vskip 20pt

\begin{center}
{{\bfseries Yude~~Liu\textsuperscript{1}\quad\quad Xinbao~~Lu\textsuperscript{2}\quad\quad Qiang~~Sun\textsuperscript{3}}\\[6pt]
\small{\textsuperscript{1}School of Mathematical Sciences, Fudan University, Shanghai, 200433, PR China\\
\textsuperscript{2}School of Mathematical Sciences, Key Laboratory of Intelligent Computing and Applications (Tongji University), Ministry of Education, Tongji University, Shanghai, 200092, PR China\\
\textsuperscript{3}Shanghai Center for Mathematical Sciences, Fudan University, Shanghai, 200438, PR China}}
\end{center}

\vskip 10pt

\footnotetext{E-mail addresses: 1.\ liuyude@fudan.edu.cn; 2.\ xinbaolu@tongji.edu.cn; 3.\ sunqiang@fudan.edu.cn.}
\footnotetext{Research of the first author is supported by China Postdoctoral Science Foundation No.~2025M773080. Research of the second author is supported by NSFC No.~12501078. Research of the third author is supported by China Postdoctoral Science Foundation No. 2025M783084 and China National Postdoctoral Program for Innovative Talents No.~BX20250071.}
\begin{center}
\begin{minipage}{15cm}
{\bfseries Abstract} We establish the Brunn--Minkowski inequality for chord integrals in $\mathbb{R}^n$, along with the equality condition, for $q\in(0,1)\cup(n+1,n+3+\eta_n]$, where $\eta_n>0$ depends only on $n$. We also prove the reverse inequality in $\mathbb{R}^2$ for $q\in(-1,0)$. As applications, we prove uniqueness of solutions to the corresponding chord Minkowski problem and derive a Urysohn inequality for fractional perimeters.


\end{minipage}
\end{center}

\vskip 25pt
\section{\bfseries Introduction}
\label{sec:statements}

The setting of this article is the $n$-dimensional Euclidean space $\mathbb{R}^n$, $n\ge2$.
A \emph{convex body} in $\mathbb{R}^n$ is a compact convex set with nonempty interior.
Denote by $\mathcal{K}^n$ the set of convex bodies in $\mathbb{R}^n$.

The classical Brunn--Minkowski inequality, which strengthens the classical
isoperimetric inequality, is a marvelous result of combining the two
basic notions of vector addition and volume.
It states that for all $K,L\in\mathcal{K}^n$,
\begin{equation}\label{intro:eq:volume-bm}
	|K+L|^{\frac1n}\ge|K|^{\frac1n}+|L|^{\frac1n},
\end{equation}
with equality if and only if $K$ and $L$ are homothetic, that is,
$L=z+\lambda K$ for some $z\in\R^n$ and $\lambda>0$.
Here $|K|$ denotes the $n$-dimensional volume of $K$.
For the connections of the Brunn--Minkowski inequality with other
geometric and analytic inequalities, see Gardner's excellent survey~\cite{Gardner02}.

Let $A(n,1)$ denote the affine Grassmannian of unoriented lines in $\R^n$.
Crofton's formula expresses the volume in \eqref{intro:eq:volume-bm} as
\[
 |K|=\int_{A(n,1)}|K\cap\ell|\dd\ell,
\]
where $|K\cap\ell|$ denotes the length of the chord $K\cap\ell$ and
$\dd\ell$ is the invariant measure on unoriented lines in $\R^n$.
Replacing $|K\cap\ell|$ by $|K\cap\ell|^q$ gives the $q$-th chord integral
\begin{equation}\label{eq:chord-definition}
 I_q(K)=\int_{|K\cap\ell|>0}|K\cap\ell|^q\dd\ell,
 \qquad q>-1.
\end{equation}
Chord integrals have long been studied in integral geometry and in the
theory of convex bodies~\cite{Hadwiger50,Santalo86,Zhang99,XS08,LXYZ}.

Cauchy's surface area formula, Crofton's volume formula, and the
Poincar\'e--Hadwiger formula (see \cite[Chapter~14]{Santalo}) give
\begin{equation}\label{intro:eq:classical-values}
 I_0(K)=\frac{\omega_{n-1}}{n\omega_n}S(K),\qquad
 I_1(K)=|K|,\qquad
 I_{n+1}(K)=\frac{n+1}{\omega_n}|K|^2,
\end{equation}
where $S(K)$ denotes the surface area of $K$ and $\omega_j$ the volume
of the Euclidean unit ball in $\R^j$.
These identities lead to the following question: for which exponents $q>-1$
does the Brunn--Minkowski inequality
\begin{equation}\label{intro:eq:chord-bm}
 I_q(K+L)^{\frac{1}{n+q-1}}
 \ge I_q(K)^{\frac{1}{n+q-1}}+I_q(L)^{\frac{1}{n+q-1}}
\end{equation}
hold for all $K,L\in\mathcal{K}^n$?
The exponent $\frac{1}{n+q-1}$ is determined by the homogeneity of $I_q$
of degree $n+q-1$.
By \eqref{intro:eq:classical-values}, the cases $q=1$ and $q=n+1$
are precisely the classical Brunn--Minkowski inequality for volume
\eqref{intro:eq:volume-bm}, while $q=0$ gives the Brunn--Minkowski
inequality for surface area.

Li~\cite{LiTalk26} announced \eqref{intro:eq:chord-bm}, together with
its equality condition, for $q\in[1,n+1]$. Li, Wan, and
Zhang~\cite[Theorem~1.1]{LWZ26} proved the inequality for
$q\in[1,n+2]\cup\{n+3\}$ in every dimension and for $q\in[1,5]$ in
$\mathbb{R}^2$, with equality precisely for homothetic bodies.

In this paper, we prove the following result.

\begin{theorem}\label{thm:chord-main}
For every $n\ge2$ there exists $\eta_n>0$ such that, for every
$0<q\le n+3+\eta_n$ and all $K,L\in\mathcal{K}^n$,
\[
 I_q(K+L)^{\frac1{n+q-1}}\ge
 I_q(K)^{\frac1{n+q-1}}+I_q(L)^{\frac1{n+q-1}},
\]
with equality if and only if $K$ and $L$ are homothetic.
\end{theorem}

For $s\in(0,1)$, the fractional perimeter of $K\in\mathcal{K}^n$ is
\[
 P_s(K)=\int_K\int_{\R^n\setminus K}|x-y|^{-n-s}\dd y\dd x.
\]
By Ludwig's slicing formula~\cite{Ludwig13} (see also Qin~\cite{Qin23}),
\[
 I_{1-s}(K)=\frac{s(1-s)}{n\omega_n}P_s(K).
\]
Thus the case $q=1-s$ of Theorem~\ref{thm:chord-main} takes the following form.

\begin{corollary}\label{cor:fractional-perimeter-main}
For every $n\ge2$, $0<s<1$, and all $K,L\in\mathcal{K}^n$,
\[
 P_s(K+L)^{\frac1{n-s}}\ge
 P_s(K)^{\frac1{n-s}}+P_s(L)^{\frac1{n-s}},
\]
with equality if and only if $K$ and $L$ are homothetic.
\end{corollary}

For $q\in(-1,0)$ and $K,L\in\mathcal{K}^2$, we prove the reverse
Brunn--Minkowski inequality
\[
 I_q(K+L)^{\frac{1}{q+1}}
 \le I_q(K)^{\frac{1}{q+1}}+I_q(L)^{\frac{1}{q+1}},
\]
with equality if and only if $K$ and $L$ are homothetic
(Theorem~\ref{thm:planar-chord-bm}). For the same range of $q$ in
dimensions $n\ge3$, both directions admit counterexamples
(Theorem~\ref{negc:thm:failures}).

Li, Wan, and Zhang~\cite[Theorem~1.6]{LWZ26} constructed counterexamples
to \eqref{intro:eq:chord-bm} for every $q>n+1+\beta_*$, where
$\beta_*\simeq2.88147$.

The reverse of \eqref{intro:eq:chord-bm} fails for every $q>0$.
Let $Q=[0,1]^n$ and, for $0<\varepsilon<1$, let
\[
 K_\varepsilon=[0,1-\varepsilon]^{n-1}\times[0,\varepsilon],
 \qquad L_\varepsilon=[0,\varepsilon]^{n-1}\times[0,1-\varepsilon].
\]
Both boxes lie in $Q$, and their chord lengths tend to zero on almost
every line as $\varepsilon\downarrow0$. By dominated convergence in
\eqref{eq:chord-definition}, $I_q(K_\varepsilon)$ and $I_q(L_\varepsilon)$
both tend to zero. Since $K_\varepsilon+L_\varepsilon=Q$, the reverse of
\eqref{intro:eq:chord-bm} would imply $I_q(Q)=0$, a contradiction.

Table~\ref{pub:tab:ranges} summarizes the results.

\begin{table}[htbp]
\centering
\small
\begin{tabularx}{\linewidth}{@{}p{0.30\linewidth}p{0.13\linewidth}XX@{}}
\toprule
Chord exponent & Dimension & Forward inequality & Reverse inequality\\
\midrule
$-1<q<0$ & $n=2$ & Fails & Holds\\
$-1<q<0$ & $n\ge3$ & Fails & Fails\\
$q=0$ & $n=2$ & Equality for all pairs & Equality for all pairs\\
$q=0$ & $n\ge3$ & Holds & Fails\\
$0<q\le n+3+\eta_n$ & $n\ge2$ & Holds & Fails\\
$n+3+\eta_n<q\le n+1+\beta_*$ & $n\ge2$ & Unresolved here & Fails\\
$q>n+1+\beta_*$ & $n\ge2$ & Fails & Fails\\
\bottomrule
\end{tabularx}
\caption{Brunn--Minkowski inequalities for chord integrals.}
\label{pub:tab:ranges}
\end{table}
In every valid nonlinear inequality in the table, equality holds if and
only if $K$ and $L$ are homothetic.

For $n\ge3$, Theorem~\ref{thm:chord-main} extends the result of Li, Wan,
and Zhang~\cite[Theorem~1.1]{LWZ26} to
$q\in(0,1)\cup(n+2,n+3)\cup(n+3,n+3+\eta_n]$.
In $\mathbb{R}^2$, it extends their range $[1,5]$ to $(0,5+\eta_2]$.
For $q\in(1,n+1)\cup(n+1,n+2)$, we also give alternative proofs.

Theorem~\ref{thm:chord-main} and the variational formula of Lutwak, Xi,
Yang, and Zhang~\cite{LXYZ} yield a first Minkowski inequality for chord
integrals. It gives uniqueness up to translation of solutions to the
chord Minkowski problem for $q\in(0,n+3+\eta_n]$, extending the result
of Li, Wan, and Zhang~\cite[Theorem~1.7(1)]{LWZ26}. For $q=1-s$ with
$s\in(0,1)$, comparison with a ball gives a Urysohn inequality for
fractional perimeters.

This article is organized as follows. In Section~\ref{sec:preliminaries},
we collect some basic facts on convex bodies and chord integrals.
We prove the planar inequalities for $0<|q|<1$ in Section~\ref{sec:planar}.
In Section~\ref{negc:sec:nonpositive}, we give counterexamples to
\eqref{intro:eq:chord-bm} and its reverse when $-1<q<0$ and $n\ge3$.
We complete the proof of Theorem~\ref{thm:chord-main} in
Sections~\ref{sec:fractional-lifting}--\ref{sec:cnd-full-interval}
and give two applications in Section~\ref{app:sec:chord-uniqueness}.

\vskip 25pt
\section{\bfseries Preliminaries}
\label{sec:preliminaries}

For quick later reference, we collect some basic facts on convex bodies. We refer to \cite{XX22} and
\cite{TXX23} for further background.

As usual, $x\cdot y$ denotes the standard inner product of $x$ and $y$ in $\R^n$.
Let $\mathbb S^{n-1}=\{u\in\R^n:|u|=1\}$ be the boundary of the
Euclidean unit ball $B^n$ in $\R^n$. The volume of $B^n$ is $\omega_n=|B^n|$.

For $u\in\mathbb S^{n-1}$, let $u^\perp=\{x\in\R^n:x\cdot u=0\}$.
For a subspace $E$ of $\R^n$, write $K|E$ for the orthogonal projection of $K$ onto $E$.

The \emph{support function} $h_K$ and \emph{radial function} $\rho_K$ of a convex body $K$ are
defined for $x\in\R^n$ by
\[
h_K(x)=\max\{x\cdot y:y\in K\},\qquad \rho_K(x)=\max\{\lambda\geq0:\lambda x\in K\}.
\]
The distance between $K,L\in\mathcal K^n$ is
$\delta_H(K,L)=\max_{u\in\mathbb S^{n-1}}|h_K(u)-h_L(u)|$.

For every nonnegative Borel function $\varphi$ on $A(n,1)$,
\[
\int_{A(n,1)}\varphi(\ell)\dd\ell
=\frac1{n\omega_n}\int_{\mathbb S^{n-1}}\int_{u^\perp}
\varphi(y+\R u)\dd y\dd u,
\]
where $\dd u$ and $\dd y$ are surface measure on $\mathbb S^{n-1}$
and Lebesgue measure on $u^\perp$, respectively.

For $K\in\mathcal K^n$, $u\in\mathbb S^{n-1}$ and $y\in u^\perp$,
the chord length of $K$ on $y+\R u$ is $\ell_K(y,u)=|K\cap(y+\R u)|$.
The function $\ell_K(\,\cdot\,,u)$ is concave and positive on
$\operatorname{int}(K|u^\perp)$. For $q>-1$, the directional chord integral is defined by
\[
Z_q(K,u)=\int_{\operatorname{int}(K|u^\perp)}\ell_K(y,u)^q\dd y.
\]
The function $Z_q(K,\cdot)$ is even. For $z\in\R^n$, $\lambda>0$,
and every rotation $T$,
\[
Z_q(z+\lambda K,u)=\lambda^{n+q-1}Z_q(K,u),\qquad Z_q(TK,Tu)=Z_q(K,u).
\]
For $q>0$, $Z_q(\,\cdot\,,u)$ is increasing under inclusion. Also $Z_0(K,u)=\mathcal H^{n-1}(K|u^\perp)$ and $Z_1(K,u)=|K|$.

Since \(\mathcal H^{n-1}(\partial(K|u^\perp))=0\), the line-measure
decomposition yields
\begin{equation}\label{eq:chord-directional}
	I_q(K)=\frac1{n\omega_n}\int_{\mathbb S^{n-1}}Z_q(K,u)\dd u.
\end{equation}

Thus \(I_q\) is rigid-motion invariant and
\[
I_q(\lambda K)=\lambda^{n+q-1}I_q(K),\qquad
K\subset L,\ q\ge0\ \Longrightarrow\ I_q(K)\le I_q(L).
\]

Write \(Z_q^{(n)}\) and \(I_q^{(n)}\) whenever the ambient dimension has
to be specified. Let \(\nu\) be a nonzero finite Borel measure on
\(\mathbb S^{n-1}\), and define
\[
J_{q,\nu}(K)=\int_{\mathbb S^{n-1}}Z_q(K,u)\,\dd\nu(u),
\qquad q>-1.
\]
Moreover,
\[
J_{q,\nu}(z+\lambda K)=\lambda^{n+q-1}J_{q,\nu}(K),\qquad
J_{q,\nu}(K)=J_{q,\nu_e}(K),
\]
where \(\nu_e(A)=\tfrac12\bigl(\nu(A)+\nu(-A)\bigr)\). Also,
\[
J_{q,\nu}=I_q
\qquad\text{when}\qquad
\dd\nu(u)=\frac{\dd u}{n\omega_n}.
\]

For \(q\ge0\), \(J_{q,\nu}(K)<\infty\). If \(rB^n\subseteq K\) and
\(\delta_H(K,L)\le\delta<r\), then
\[
\left(1-\frac{\delta}{r}\right)K
\subseteq L\subseteq
\left(1+\frac{\delta}{r}\right)K.
\]
If \(G\) is finite, nonnegative, inclusion monotone, and
\(G(\lambda K)=\lambda^\alpha G(K)\), \(\alpha>0\), then
\begin{equation}\label{eq:monotone-continuity}
	\left(1-\frac{\delta}{r}\right)^\alpha G(K)
	\le G(L)\le
	\left(1+\frac{\delta}{r}\right)^\alpha G(K).
\end{equation}
If \(G\) is translation invariant, \eqref{eq:monotone-continuity} gives
Hausdorff continuity. For \(K_j\to K\) and \(L_j\to L\),
the uniform convergence
\[
\sup_{0\le t\le1}
\left|
G\bigl((1-t)K_j+tL_j\bigr)
-
G\bigl((1-t)K+tL\bigr)
\right|
\longrightarrow0.
\]

For \(q>0\), apply \eqref{eq:monotone-continuity} to
\(G=Z_q(\,\cdot\,,u)\). The constants are independent of \(u\).
Hence \(J_{q,\nu}\) is continuous. For \(-1<q<0\), use the next lemma.

\begin{lemma}\label{lem:continuity}
	For $-1<q<0$, the function $Z_q(K,u)$ is finite and positive, and is jointly continuous in $(K,u)\in\mathcal K^n\times\mathbb S^{n-1}$. Consequently $J_{q,\nu}$ is finite, positive, and continuous for every nonzero finite Borel measure $\nu$ on $\mathbb S^{n-1}$. 
\end{lemma}
\begin{proof}
	Fix \(u\in\mathbb S^{n-1}\) and let \(K_j\to K\). Translate all bodies
	once and choose \(0<b_1<b_2\) so that, for large \(j\),
	\[
	b_1B^n\subset K_j,K\subset b_2B^n.
	\]
	Put \(D_j=K_j|u^\perp\) and \(D=K|u^\perp\). Concavity of the
	section length, together with \(\ell_{K_j}(0,u)\ge2b_1\), gives
	\[\ell_{K_j}(t\rho_{D_j}(v)v,u)\ge2b_1(1-t),\quad \qquad~\text{for}~v\in\mathbb S^{n-1}\cap u^\perp,\  t\in [0,1].\]
	Then
	\[
	\begin{aligned}
		0<Z_q(K_j,u)
		&=\int_{\mathbb S^{n-1}\cap u^\perp}\int_0^1
		\ell_{K_j}\bigl(t\rho_{D_j}(v)v,u\bigr)^q
		\rho_{D_j}(v)^{n-1}t^{n-2}\,dt\,dv\\
		&\leq b_2^{n-1}(2b_1)^q
		\int_{\mathbb S^{n-1}\cap u^\perp}\int_0^1
		t^{n-2}(1-t)^q\,dt\,dv,
	\end{aligned}
	\]
	The bound is finite for \(q>-1\).
	
	Choose \(\varepsilon_j\downarrow0\) so that
	\[
	(1-\varepsilon_j)K\subseteq K_j\subseteq(1+\varepsilon_j)K,
	\qquad
	(1-\varepsilon_j)\rho_D(v)\leq\rho_{D_j}(v)
	\leq(1+\varepsilon_j)\rho_D(v).
	\]
	Thus,
	\[
	\|\rho_{D_j}-\rho_D\|_{\infty}\longrightarrow0.
	\]
	For fixed \(0\le t<1\), and large \(j\),
	\[
	(1-\varepsilon_j)
	\ell_K\left(
	\frac{t\rho_{D_j}(v)}{1-\varepsilon_j}v,u
	\right)
	\le
	\ell_{K_j}\bigl(t\rho_{D_j}(v)v,u\bigr)
	\le
	(1+\varepsilon_j)
	\ell_K\left(
	\frac{t\rho_{D_j}(v)}{1+\varepsilon_j}v,u
	\right).
	\]
	Since
	\[
	\frac{t\rho_{D_j}(v)}{1\pm\varepsilon_j}v
	\longrightarrow t\rho_D(v)v\in\operatorname{relint}D,
	\]
	continuity of chord length gives
	\[
	\ell_{K_j}\bigl(t\rho_{D_j}(v)v,u\bigr)
	\longrightarrow \ell_K\bigl(t\rho_D(v)v,u\bigr).
	\]
	Dominated convergence gives
	\[
	\begin{aligned}
		Z_q(K_j,u)
		&=\int_{\mathbb S^{n-1}\cap u^\perp}\int_0^1
		\ell_{K_j}\bigl(t\rho_{D_j}(v)v,u\bigr)^q
		\rho_{D_j}(v)^{n-1}t^{n-2}\,dt\,dv\\
		&\longrightarrow
		\int_{\mathbb S^{n-1}\cap u^\perp}\int_0^1
		\ell_K\bigl(t\rho_D(v)v,u\bigr)^q
		\rho_D(v)^{n-1}t^{n-2}\,dt\,dv
		=Z_q(K,u).
	\end{aligned}
	\]
	
	For \(u_j\to u\), choose rotations \(T_j\to\mathrm{Id}\) with
	\(T_ju=u_j\). Then \(T_j^{-1}K_j\to K\), and rotation covariance gives
	\[
	Z_q(K_j,u_j)
	=Z_q(T_j^{-1}K_j,T_j^{-1}u_j)
	=Z_q(T_j^{-1}K_j,u)
	\longrightarrow Z_q(K,u).
	\]
	Thus, by compactness of the sphere,
	\[
	\sup_{u\in\mathbb S^{n-1}}
	\left|Z_q(K_j,u)-Z_q(K,u)\right|\longrightarrow0.
	\]
	Hence
	\[
	\left|J_{q,\nu}(K_j)-J_{q,\nu}(K)\right|
	\leq\nu(\mathbb S^{n-1})
	\sup_{u\in\mathbb S^{n-1}}
	\left|Z_q(K_j,u)-Z_q(K,u)\right|
	\longrightarrow0.
	\]
	This complete the proof.
\end{proof}

Write $K^2=K\times K$. If $q>1$ and $\dd\nu(u)=f(u)\dd u$,
where $f$ is nonnegative and integrable on $\mathbb S^{n-1}$, then
for every $K\in\mathcal K^n$,
\begin{equation}\label{eq:weighted-pair-energy}
	J_{q,\nu}(K)=q(q-1)\int_{K^2}
	f\Bigl(\frac{y-x}{|y-x|}\Bigr)|x-y|^{q-n-1}\dd x\dd y.
\end{equation}
To prove this, let $g_K(z)=|K\cap(K+z)|$.
Since every section parallel to $u$ is an interval, for $r\ge0$,
\begin{align}
	g_K(ru)&=\int_{u^\perp}(\ell_K(y,u)-r)_+\dd y,
	\label{eq:covariogram-slicing}\\
	|K|-g_K(ru)&=\int_{u^\perp}\min\{r,\ell_K(y,u)\}\dd y.
	\label{eq:covariogram-deficit-slicing}
\end{align}
Equation \eqref{eq:covariogram-slicing} and Tonelli's theorem give,
for each $u\in\mathbb S^{n-1}$,
\[
\int_0^\infty r^{q-2}g_K(ru)\dd r=\int_{u^\perp}\int_0^{\ell_K(y,u)}
r^{q-2}(\ell_K(y,u)-r)\dd r\dd y=\frac{Z_q(K,u)}{q(q-1)}.
\]
For any nonnegative Borel function $\varphi$ on $\R^n$, the substitution
$z=y-x$ gives
\begin{equation}\label{eq:covariogram-change-variables}
	\int_{K^2}\varphi(y-x)\dd x\dd y=\int_{\R^n}\varphi(z)|K\cap(K-z)|\dd z=\int_{\R^n}\varphi(z)g_K(z)\dd z,
\end{equation}
where the last equality follows by translating $K\cap(K-z)$ by $z$.
The identity also holds for signed Borel functions when the integrals converge absolutely.
Using the directional identity above and then polar coordinates $z=ru$,
we obtain
\[
\begin{aligned}
	J_{q,\nu}(K)
	&=\int_{\mathbb S^{n-1}}Z_q(K,u)f(u)\dd u=q(q-1)\int_{\mathbb S^{n-1}}f(u)
	\int_0^\infty r^{q-2}g_K(ru)\dd r\dd u\\
	&=q(q-1)\int_{\R^n}
	f\Bigl(\frac z{|z|}\Bigr)|z|^{q-n-1}g_K(z)\dd z.
\end{aligned}
\]
This, together with \eqref{eq:covariogram-change-variables}, yields
\eqref{eq:weighted-pair-energy}.

Taking $f=\frac{1}{n\omega_n}$ in \eqref{eq:weighted-pair-energy} gives
\begin{equation}\label{intro:eq:chord-energy}
	I_q(K)=\frac{q(q-1)}{n\omega_n}
	\int_{K^2}|x-y|^{q-n-1}\dd x\dd y,\qquad q>1.
\end{equation}

For $K\in\mathcal K^n$, define
\[
\begin{aligned}
 E_a(K)&=\int_{K^2}|x-y|^{-a}\dd x\dd y, &&0\le a<n,\\
 \mathcal E_r^{(n)}(K)&=\int_{K^2}|x-y|^r\dd x\dd y, &&r>0.
\end{aligned}
\]
Their homogeneity and monotonicity under inclusion give 
continuity by \eqref{eq:monotone-continuity}.

\begin{lemma}\label{lem:homogeneous-bm}
	Let $F:\mathcal K^n\to(0,\infty)$ be homogeneous of degree one.
	Then the following statements are equivalent.
	\begin{enumerate}
		\item[(i)] $F(K+L)\ge F(K)+F(L)$ for all $K,L\in\mathcal K^n$.
		\item[(ii)] $F((1-t)K+tL)\ge(1-t)F(K)+tF(L)$ for all $K,L\in\mathcal K^n$
		and $t\in[0,1]$.
		\item[(iii)] $F((1-t)K+tL)\ge1$ for all $K,L\in\mathcal K^n$ with $F(K)=F(L)=1$
		and all $t\in[0,1]$.
	\end{enumerate}
	If these inequalities hold, (i) and (ii) have the same equality cases
	for every $t\in(0,1)$.
	The same assertions hold with all inequalities reversed.
\end{lemma}

\begin{proof}
	Write $K_t=(1-t)K+tL$. For $0<t<1$, (i) and homogeneity give
	\[
	F(K_t)\ge F((1-t)K)+F(tL)=(1-t)F(K)+tF(L).
	\]
	At $t=0,1$, (ii) holds with equality. Thus (i)$\Rightarrow$(ii), while (ii)$\Rightarrow$(iii) is immediate.
	Conversely, set
	\[
	a=F(K),\qquad b=F(L),\qquad A=(1-t)a+tb.
	\]
	Since $F\bigl(\frac{K}{a}\bigr)=F\bigl(\frac{L}{b}\bigr)=1$, (iii) and homogeneity give
	\[
	F(K_t)=A F\left(\frac{(1-t)a}{A}\frac Ka+\frac{tb}{A}\frac Lb\right)
	\ge A=(1-t)F(K)+tF(L).
	\]
	Thus (iii)$\Rightarrow$(ii). Taking $t=\frac{1}{2}$ in (ii) yields
	\[
	F(K+L)=2F\left(\frac{K+L}{2}\right)\ge F(K)+F(L),
	\]
	and hence (ii)$\Rightarrow$(i).
	
	Assume first that the inequalities have the displayed direction. The
	function
	\[
	G(t)=F(K_t)-(1-t)F(K)-tF(L)
	\]
	is then concave, nonnegative, and zero at both endpoints. If it vanishes
	at one \(t_0\in(0,1)\), it vanishes on the whole interval. Thus an
	interior equality in (ii) is equivalent to equality at \(t=\frac{1}{2}\), and
	hence to equality in (i). The argument with all signs reversed is the
	same, with convexity in place of concavity.
\end{proof}

\vskip 25pt
\section{\bfseries Planar Brunn--Minkowski inequalities}
\label{sec:planar}
\label{simp:sec:planar-shared}

We prove the Brunn--Minkowski inequality for chord integrals in $\mathbb{R}^2$
when $0<q<1$ and its reverse when $-1<q<0$, including their equality cases.
We first establish the following more general inequalities for $J_{q,\nu}$.

\begin{theorem}
	\label{negp:thm:reverse}
	Let $\nu$ be a nonzero finite Borel measure on $\mathbb{S}^1$.
	If $0<q<1$, then for $K,L\in\mathcal{K}^2$,
	\[
	J_{q,\nu}(K+L)^{\frac1{q+1}}\ge
	J_{q,\nu}(K)^{\frac1{q+1}}+J_{q,\nu}(L)^{\frac1{q+1}}.
	\]
	If $-1<q<0$, then for $K,L\in\mathcal{K}^2$,
	\[
	J_{q,\nu}(K+L)^{\frac1{q+1}}\le
	J_{q,\nu}(K)^{\frac1{q+1}}+J_{q,\nu}(L)^{\frac1{q+1}}.
	\]
	In each case, equality holds if and only if $K$ and $L$ are homothetic.
\end{theorem}
\begin{proof}
	Put $K_t=(1-t)K+tL$ and
	$F(t)=J_{q,\nu}(K_t)^{\frac1{q+1}}$.
	By Lemma~\ref{lem:homogeneous-bm}, it suffices to prove
	\[
	q\bigl[F(t)-(1-t)F(0)-tF(1)\bigr]\ge0,\qquad 0<t<1,
	\]
	with equality if and only if $K$ and $L$ are homothetic. Set
	\[
	W(t)=|K_t|^{\frac12},\qquad
	c=\min_{0\le t\le1}\frac{F(t)}{W(t)}>0.
	\]
	To establish the equality condition as well, we prove the stronger estimate
	\begin{equation}\label{simp:eq:finite-area}
		\begin{aligned}
			q\bigl[F(t)-(1-t)F(0)-tF(1)\bigr]
			&\ge q^2c\bigl[W(t)-(1-t)W(0)-tW(1)\bigr].
		\end{aligned}
	\end{equation}
	By the Brunn--Minkowski inequality for volume
	\cite[Theorem~4.1 and pp.~363--364]{Gardner02}, the right-hand side is
	nonnegative and vanishes if and only if $K$ and $L$ are homothetic.
	
	\textbf{Step~1.} Prove a variational formula for chord integrals.\par\nopagebreak
	
	Since $J_{q,\nu}$ depends only on the even part of $\nu$, we may assume
	that $\nu$ is even. We first consider measures of the form
	$\dd\nu(u)=f(u)\,\dd u$,
	where $f$ is a nonzero smooth nonnegative even function on $\mathbb{S}^1$.
	
	Fix distinct unit vectors $u_1,\ldots,u_m$ in counterclockwise order,
	not contained in any closed semicircle.
	For $h=(h_1,\ldots,h_m)\in\mathbb R^m$, write
	\[
	P(h)=\bigcap_{i=1}^m\{x:x\cdot u_i\le h_i\}.
	\]
	We consider the open set of parameters $h$ for which $P(h)$ has $m$ edges.
	Denote the edge with outer normal $u_i$ by $E_i(h)$,
	with indices understood~cyclically.
	
	The common vertex of $E_i(h)$ and $E_{i+1}(h)$ is determined by
	$x\cdot u_i=h_i$ and $x\cdot u_{i+1}=h_{i+1}$.
	Since $u_i$ and $u_{i+1}$ are fixed and linearly independent,
	this vertex depends linearly on $h$.
	The edge lengths are then also linear in $h$, since the edge directions are~fixed.
	
	We prove that, for $q>-1$, the map $h\mapsto J_{q,\nu}(P(h))$ is $C^2$
	and that, for nonadjacent edges $E_i(h)$ and $E_j(h)$,
	\begin{equation}\label{simp:eq:weighted-offdiag}
		\frac{\partial^2 J_{q,\nu}(P(h))}{\partial h_i\,\partial h_j}
		=2q(q-1)\int_{E_i(h)}\int_{E_j(h)}
		f\Bigl(\frac{y-x}{|y-x|}\Bigr)|x-y|^{q-3}
		\dd y\dd x.
	\end{equation}
	In integrals over line segments, $\dd x$ and $\dd y$ denote arc-length elements.
	
	First, let $\psi:[0,\infty)\to[0,\infty)$ be any smooth function
	that vanishes near zero.
	Set
	\[
	J_\psi(P(h))=\int_{\mathbb S^1}\int_{u^\perp}
	\psi\bigl(\ell_{P(h)}(y,u)\bigr)\dd y\,f(u)\dd u.
	\]
	Here $\ell_{P(h)}(y,u)=|P(h)\cap(y+\R u)|$,
	as in Section~\ref{sec:preliminaries}.
	
	To differentiate $J_\psi(P(h))$, we first express it as an integral over $P(h)^2$.
	Since $\psi(0)=\psi'(0)=0$, integration by parts gives
	\[
	\int_0^\ell(\ell-r)\psi''(r)\dd r
	=\bigl[(\ell-r)\psi'(r)\bigr]_0^\ell
	+\psi(\ell)-\psi(0)=-\ell\psi'(0)+\psi(\ell)-\psi(0)
	=\psi(\ell).
	\]
	Substituting the preceding identity into the definition of $J_\psi(P(h))$,
	interchanging the integrals, using \eqref{eq:covariogram-slicing},  polar coordinates and
	\eqref{eq:covariogram-change-variables}, we have
	\[
	J_\psi(P(h))
	=\int_{\mathbb S^1}f(u)\int_0^\infty
	g_{P(h)}(ru)\psi''(r)\dd r\dd u=\int_{P(h)^2}
	f\Bigl(\frac{y-x}{|y-x|}\Bigr)
	\frac{\psi''(|x-y|)}{|x-y|}\dd x\dd y.
	\]
	
	Define
	\[
	\Phi_\psi(x,y)=
	\begin{cases}
		\displaystyle f\Bigl(\frac{y-x}{|y-x|}\Bigr)
		\frac{\psi''(|x-y|)}{|x-y|},
		&x\neq y,\\[2mm]
		0,&x=y,
	\end{cases}
	\]
	
	Since $f,\psi$ is smooth and $\psi$ vanishes near zero, $\Phi$ is smooth.
	
	In a sufficiently small neighborhood of $h$, triangulate each $P(h)$ by joining a fixed common
	interior point to its vertices. Let $T_{i,h}$, $1\le i\le m$, be affine maps
	from a fixed reference triangle $\Delta$ onto these triangles. Then
	\[ J_\psi(P(h))=\sum_{i,j=1}^{m}|\det DT_{i,h}|\,|\det DT_{j,h}|\int_{\Delta\times\Delta}
	\Phi\bigl(T_{i,h}(\xi),T_{j,h}(\eta)\bigr)\dd\xi\dd\eta.
	\]
	The polygon vertices depend linearly on $h$, so the maps $T_{i,h}$ and their
	Jacobians depend smoothly on $h$. Since $\Delta\times\Delta$ is fixed and
	compact, we may differentiate twice under the integral sign, proving that
	$J_\psi(P(h))$ is $C^2$ in $h$.
	
	Fix $i\in\{1,\ldots,m\}$. Keeping $h_k$ fixed for $k\ne i$, set
	$Q_i=\bigcap_{k\ne i}\{x:x\cdot u_k\le h_k\}$. Then
	$P(h)=Q_i\cap\{x:x\cdot u_i\le h_i\}$.
	Let $\varphi$ be a continuous function on $\mathbb R^2$ independent of $h$.
	Since the endpoints of $Q_i\cap\{x:x\cdot u_i=r\}$ vary continuously
	with $r$ near $h_i$, the integral of $\varphi$ over this segment is
	continuous. Thus,
	\[
	\frac{\partial}{\partial h_i}\int_{P(h)}\varphi(x)\dd x=\frac{\partial}{\partial h_i}
	\int_{-\infty}^{h_i}
	\int_{Q_i\cap\{x:x\cdot u_i=r\}}\!\!\!\!\!\!\!\!\!
	\varphi(x)\dd x\dd r=\int_{Q_i\cap\{x:x\cdot u_i=h_i\}}
	\varphi(x)\dd x=\int_{E_i(h)}\varphi(x)\dd x.
	\]
	Applying this identity to both integration domains and using the
	evenness of $f$, we obtain, with $r=|x-y|$,
	\[
	\begin{aligned}
		\frac{\partial}{\partial h_i}J_\psi(P(h))
		&=\int_{E_i(h)}\int_{P(h)}
		f\Bigl(\frac{y-x}{r}\Bigr)\frac{\psi''(r)}r
		\dd y\dd x\\
		&\quad+\int_{P(h)}\int_{E_i(h)}
		f\Bigl(\frac{y-x}{r}\Bigr)\frac{\psi''(r)}r
		\dd y\dd x\\
		&=2\int_{E_i(h)}\int_{P(h)}
		f\Bigl(\frac{y-x}{r}\Bigr)\frac{\psi''(r)}r
		\dd y\dd x.
	\end{aligned}
	\]
	For nonadjacent edges $E_i(h)$ and $E_j(h)$, the endpoints of $E_i(h)$
	depend only on $h_{i-1}$, $h_i$ and $h_{i+1}$, so $E_i(h)$ remains fixed as $h_j$ varies. Then
	\begin{equation}\label{eq:smooth-chord-offdiag}
		\frac{\partial^2J_\psi(P(h))}{\partial h_i\,\partial h_j}
		=2\int_{E_i(h)}\int_{E_j(h)}
		f\Bigl(\frac{y-x}{r}\Bigr)\frac{\psi''(r)}r
		\dd y\dd x.
	\end{equation}
	
	Now fix $q>-1$. We will choose $\psi$ to agree with $r^q$ away from zero
	and show that $J_{q,\nu}(P(h))-J_\psi(P(h))$ is locally constant in $h$.
	
	Denote $\ell_k(h)$ as the length of $E_k(h)$, $1\le k\le m$. For every adjacent pair of edges $E_i(h), E_j(h)$, let \(x_0(h)=E_i(h)\cap E_j(h)\) and
	write \(x=x_0(h)+\xi v_i\) and \(y=x_0(h)+\eta v_j\), where
	\(\xi,\eta\ge0\)  $v_i,v_j$ are the  unit vectors pointing along $E_i(h), E_j(h)$ away from $x_0(h)$.  There exists \(c_0>0\), independent of \(h\),
	such that \(|x-y|\ge c_0(\xi+\eta)\). 
	
	The edge lengths and the distances between nonadjacent edges are positive
	and continuous in $h$. Around any point of this open set, fix a neighborhood
	on which these finitely many quantities have a common positive lower bound.
	Choose $\varepsilon>0$ so that, throughout this neighborhood,
	\[
	\begin{aligned}
		d(E_i(h),E_j(h))&>2\varepsilon,\quad\text{for all nonadjacent pairs},\\
		c_0\ell_k(h)&>2\varepsilon,\quad\text{for }k=1,\ldots,m.
	\end{aligned}
	\]
	Choose $\psi$ so that $0\le\psi(r)\le r^q$ for $r>0$,
	$\psi(r)=0$ for $r\le\varepsilon$, and $\psi(r)=r^q$ for $r\ge2\varepsilon$.
	
	This choice gives
	$0\le J_\psi(P(h))\le J_{q,\nu}(P(h))<\infty$.
	With $r=|x-y|$, Z\"ahle's boundary Blaschke--Petkantschin formula
	\cite[formula~(6), p.~176]{Moseeva22}~gives
	\begin{equation}\label{eq:boundary-chord}
		J_{q,\nu}(P(h))-J_\psi(P(h))
		=\sum_{i\ne j}\int_{E_i(h)}\int_{E_j(h)}
		f\Bigl(\frac{y-x}{r}\Bigr)\frac{r^q-\psi(r)}{r^3}[u_i\cdot(x-y)][u_j\cdot(y-x)]
		\dd y\dd x.
	\end{equation}
	For $x\in E_i(h)$ and $y\in E_j(h)$, we have
	\[
	u_i\cdot(x-y)=h_i-u_i\cdot y\ge0,\qquad
	u_j\cdot(y-x)=h_j-u_j\cdot x\ge0.
	\]
	Thus the integrals on the right-hand side of \eqref{eq:boundary-chord}
	are nonnegative. Since the left-hand side is finite, each integral is finite.
	
	For nonadjacent edges, $r>2\varepsilon$ by the choice of $\varepsilon$,
	so $\psi(r)=r^q$ and the integrals over these pairs in
	\eqref{eq:boundary-chord} vanish.
	Thus only adjacent pairs contribute to the sum.
	
	For an adjacent pair, the integrand in \eqref{eq:boundary-chord}
	depends on $x,y$ only through $y-x$.
	Under the parametrization above, $y-x=\eta v_j-\xi v_i$,
	so the integrand is independent of $h$ in these coordinates.
	Since $r\ge c_0(\xi+\eta)$ and $r^q-\psi(r)=0$ for $r\ge2\varepsilon$,
	the integrand vanishes
	outside the fixed triangle
	$\{\xi,\eta\ge0:\xi+\eta<2\varepsilon/c_0\}$.
	By the choice of $\varepsilon$, this triangle is contained in
	$[0,\ell_i(h)]\times[0,\ell_j(h)]$.
	Each integral over adjacent edges can therefore be taken over this fixed triangle,
	with an integrand independent of $h$.
	It follows that $J_{q,\nu}(P(h))-J_\psi(P(h))$ is locally~constant.
	
	Consequently, $J_{q,\nu}(P(h))$ is $C^2$, and its mixed second derivatives
	agree with those of $J_\psi(P(h))$ on this neighborhood.
	For nonadjacent edges, $r=|x-y|>2\varepsilon$, so
	$\psi''(r)/r=q(q-1)r^{q-3}$.
	Substituting this into \eqref{eq:smooth-chord-offdiag}
	proves \eqref{simp:eq:weighted-offdiag}.
	
	\textbf{Step~2.} Compare the area Hessian.\par\nopagebreak
	
	Let $\delta_i\in(0,\pi)$ be the angle from \(u_i\) to \(u_{i+1}\),
	with cyclic indices. Denote by \(U\in\mathbb R^{m\times2}\) the matrix
	whose \(i\)-th row is \(u_i^\top\). The area Hessian $V=D^2|P(h)|$ has entries
	\[
	V_{i,i+1}=V_{i+1,i}=\csc\delta_i,\qquad
	V_{ii}=-\cot\delta_{i-1}-\cot\delta_i,
	\]
	with other entries zero.
	
	Suppose $A$ is a real symmetric $m\times m$ matrix such that
	\[
	AU=0,\qquad A_{ij}\le0\quad\text{for all nonadjacent distinct }i,j.
	\]
	Then $A\le\gamma V$ for some $\gamma\in\mathbb R$. In the following,
	we prove this by subtracting a negative semidefinite matrix that
	cancels the nonadjacent entries of $A$.
	
	For $0<\theta<\pi$, define
	\[
	Q_\theta(a,b)=\frac{2ab-\cos\theta(a^2+b^2)}{\sin\theta}.
	\]
	Fix an unordered nonadjacent pair $\{i,j\}$. We define the associated matrix $S^{ij}\in\mathbb{R}^{m\times m}$ separately in the following two cases.
	
	{\bf{Case 1.}} If $u_j\neq-u_i$, then shorter arc from $u_i$ to $u_j$ has length $\Theta<\pi$. Wirite the successive nodes of the shorter arc as
	$\alpha_0<\cdots<\alpha_k$, $k\ge2$.
	Define the quadratic form as following
	\[
	w^{\top}S^{ij}w=\sum_{r=0}^{k-1}Q_{\alpha_{r+1}-\alpha_r}(w_r,w_{r+1})
	-Q_\Theta(w_0,w_k),
	Q_\delta(a,b)=\frac{2ab-\cos\delta(a^2+b^2)}{\sin\delta}.
	\]
	Here $w\in\mathbb R^m$ is arbitrary, and $w_r$ denotes its component
	corresponding to the node $\alpha_r$; and 
	\(
	Q_\theta(a,b)=\frac{2ab-\cos\theta(a^2+b^2)}{\sin\theta}.
	\) For $a,b,c\in\mathbb{R}$ and $\alpha,\beta>0$ with $\alpha+\beta<\pi$, we have
	\[Q_\alpha(a,b)+Q_\beta(b,c)-Q_{\alpha+\beta}(a,c)=-\frac{\sin(\alpha+\beta)}{\sin\alpha\sin\beta}
	\left(b-\frac{\sin\beta\,a+\sin\alpha\,c}{\sin(\alpha+\beta)}\right)^2.
	\]
	Applying this identity successively to merge neighboring intervals gives
	$S^{ij}\le0$, since their total length is $\Theta<\pi$.
	For sine and cosine data, every squared term vanishes, hence $S^{ij}U=0$.
	Its only nonadjacent mixed entry is $S^{ij}_{ij}=-\csc\Theta$.
	
	{\bf{Case 2.}} If \(u_j=-u_i\), define
	\(S^{ij}=-(e_i+e_j)(e_i+e_j)^\top\). Then \(S^{ij}\le0\), and
	\(S^{ij}U=-(e_i+e_j)(u_i+u_j)^\top=0\). Moreover, all entries of
	\(S^{ij}\) vanish except the diagonal entries
	\(S^{ij}_{ii}=S^{ij}_{jj}=-1\) and the endpoint entries
	\(S^{ij}_{ij}=S^{ij}_{ji}=-1\).
	
	For every unordered nonadjacent pair with \(A_{ij}<0\), multiply
	\(S^{ij}\) by \(-A_{ij}\sin\Theta\) when \(\Theta<\pi\), and by
	\(-A_{ij}\) in the antipodal case. Let \(S\) be the sum of the resulting
	matrices. Since \(A_{ij}\le0\), we have \(S\le0\), \(SU=0\), and
	\(S_{ij}=A_{ij}\) for every nonadjacent pair.
	
	Set \(B=A-S\). Then \(A\le B\), \(BU=0\), and \(B_{ij}=0\) for all
	nonadjacent pairs. Hence \(B\) is cyclic tridiagonal. Since
	\(u_{i+1}=\cos\delta_i\,u_i+\sin\delta_i\,v_i\) and
	\(u_{i-1}=\cos\delta_{i-1}\,u_i-\sin\delta_{i-1}\,v_i\), where \(v_i\) is the counterclockwise unit tangent to \(u_i\), the equation
	\(BU=0\) gives
	\(B_{i,i+1}\sin\delta_i=B_{i,i-1}\sin\delta_{i-1}=B_{i-1,i}\sin\delta_{i-1}\). Thus there is
	\(\gamma\in\mathbb R\) such that
	\(B_{i,i+1}=\gamma\csc\delta_i\). Comparing the \(u_i\)-components in
	\(BU=0\) yields
	\(B_{ii}=-\gamma(\cot\delta_{i-1}+\cot\delta_i)\). Therefore
	\(B=\gamma V\), and hence
	\(A\le\gamma V\).
	
	\textbf{Step~3.} Prove \eqref{simp:eq:finite-area} for polygons.\par\nopagebreak
	
	Let $K,L$ be polygons, and list the outer normals of
	both polygons as $u_1,\ldots,u_m$ in counterclockwise
	order. Write
	\[
	K_t=P(h(t)),\qquad
	h_i(t)=(1-t)h_K(u_i)+th_L(u_i).
	\]
	Then every edge
	of $K_t$ has positive length for $0<t<1$.
	
	Set $J=J_{q,\nu}(P(h))$, $F=J^{\frac1{q+1}}$,
	$W=|P(h)|^{\frac12}$ and $A=q D_h^2J$.
	Translation invariance and homogeneity give $AU=0$ and
	$h^{\top}Ah=(q+1)q^2J$. Moreover, \eqref{simp:eq:weighted-offdiag}
	gives $A_{ij}\le0$ for nonadjacent distinct $i,j$.
	Step~2 and $h^{\top}Vh=2W^2$ imply
	\[
	A\le\gamma V,\qquad
	\gamma\ge\frac{(q+1)q^2J}{2W^2}>0.
	\]

	Set $z=h'-(W'/W)h$, where primes denote derivatives in $t$.
	Homogeneity gives $D_h^2F[h,\cdot]=0$, hence $F''=D_h^2F[z,z]$.
	We have $h^{\top}Vz=0$ and
	$z^{\top}Vz=2WW''\le0$ by the Brunn--Minkowski inequality for volume. $A\le\gamma V$, $V\le0$ and the lower bound for $\gamma$ give
	\[
	q F''=\frac{F}{(q+1)J}z^{\top}Az
	-\frac{q^2F}{(q+1)^2J^2}\bigl(D_hJ[z]\bigr)^2\le\frac{\gamma F}{(q+1)J}z^{\top}Vz
	\le q^2\frac{F}{W}W''
	\]
	Since $W''\le0$ and $F/W\ge c$, we obtain
	$(q F-q^2cW)''\le0$. This, combining with the continuity of $F,W$, gives
	\eqref{simp:eq:finite-area}.
	
	\textbf{Step~4.} Approximate and prove the equality.\par\nopagebreak
	
	For arbitrary $K,L\in\mathcal{K}^2$, choose polygons $K_j\to K$, $L_j\to L$, and set
	$K_{j,t}=(1-t)K_j+tL_j$.
	Since \(J_{q,\nu}\) depends only on the even part of \(\nu\), w.l.o.g., we can assume that $\nu$ is even.
	Choose measures \(\nu_j\) of mass \(\nu(\mathbb S^1)\) with smooth
	nonnegative even densities such that \(\nu_j\rightharpoonup\nu\).
	Set $F_j(t)=J_{q,\nu_j}(K_{j,t})^{\frac1{q+1}}$
	and $W_j(t)=|K_{j,t}|^{\frac12}$.
	By \eqref{eq:monotone-continuity} for $q>0$ and
	Lemma~\ref{lem:continuity} for $q<0$, $Z_q$ is jointly continuous.
	Thus
	\[
	\sup_{(t,u)\in[0,1]\times\mathbb S^1}
	|Z_q(K_{j,t},u)-Z_q(K_t,u)|\longrightarrow0.
	\]
	Weak convergence of $\nu_j$ is uniform on the compact family
	$\{Z_q(K_t,\cdot):0\le t\le1\}\subset C(\mathbb S^1)$. Together with
	area continuity, this gives
	\[
	\|F_j-F\|_{C[0,1]}+\|W_j-W\|_{C[0,1]}\longrightarrow0,
	\qquad
	\min_{[0,1]}\frac{F_j}{W_j}\longrightarrow
	\min_{[0,1]}\frac{F}{W}=c>0.
	\]
	Therefore \eqref{simp:eq:finite-area} passes to the limit.
	
	Equality in the claimed inequality forces equality in the
	Brunn--Minkowski inequality for volume by \eqref{simp:eq:finite-area},
	hence homothety.
	The converse follows from translation invariance and homogeneity.
\end{proof}
\begin{theorem}\label{thm:planar-chord-bm}
	Let $K,L\in\mathcal K^2$. If $0<q<1$, then
	\[
	I_q(K+L)^{\frac1{q+1}}\ge
	I_q(K)^{\frac1{q+1}}+I_q(L)^{\frac1{q+1}}.
	\]
	If $-1<q<0$, then
	\[
	I_q(K+L)^{\frac1{q+1}}\le
	I_q(K)^{\frac1{q+1}}+I_q(L)^{\frac1{q+1}}.
	\]
	In each case, equality holds if and only if $K$ and $L$ are homothetic.
\end{theorem}
\begin{proof}
	Take $\dd\nu=\dd u/(2\pi)$ in Theorem~\ref{negp:thm:reverse},
	so that $J_{q,\nu}=I_q$.
\end{proof}

\vskip 25pt
\section{\bfseries \texorpdfstring{$-1<q<0$}{-1 < q < 0}}
\label{negc:sec:nonpositive}

\begin{theorem}
	\label{negc:thm:failures}
	Let $-1<q<0$ and $n\ge 3$.
	There exist origin-symmetric $K_1,L_1,K_2,L_2\in\mathcal{K}^n$  with
	$C^\infty$ boundaries and positive Gaussian curvature such that
	\[
	I_q(K_1+L_1)^{\frac1{n+q-1}}
	<I_q(K_1)^{\frac1{n+q-1}}+I_q(L_1)^{\frac1{n+q-1}},
	\]
	and
	\[
	I_q(K_2+L_2)^{\frac1{n+q-1}}
	>I_q(K_2)^{\frac1{n+q-1}}+I_q(L_2)^{\frac1{n+q-1}}.
	\]
\end{theorem}
\begin{proof}
	For $T\ge1$, let $\varepsilon=T^{-2}$ and define the ellipsoids
	\[
	E_{\varepsilon,j}
	=(\varepsilon^2\mathrm{Id}+e_j\otimes e_j)^{\frac{1}{2}}B^n,
	\qquad j=1,\ldots,n.
	\]
	For $k\in\{1,2\}$ and $\lambda\in\{1,\frac{3}{2},2\}$, let
	\[
	A_{k,T}(\lambda)
	=\sum_{j=1}^k E_{\varepsilon,j}
	+\lambda T\sum_{j=k+1}^n E_{\varepsilon,j},
	\qquad
	K_k=A_{k,T}(1),\qquad L_k=A_{k,T}(2).
	\]
	Then $K_1,L_1,K_2,L_2$
	are origin-symmetric, with $C^\infty$ boundaries and positive Gaussian
	curvature, and
	\[
	K_k+L_k=2A_{k,T}\bigl(\frac{3}{2}\bigr).
	\]
	
	For fixed $k$ and $\lambda$, write
	\[
	D=\operatorname{diag}(\mathrm{Id}_k,\lambda T\mathrm{Id}_{n-k}),
	\qquad M=D^{-1}A_{k,T}(\lambda),
	\qquad C=[-1,1]^n.
	\]
	For $x\in\mathbb R^n$,
	\[
	0\le h_{E_{\varepsilon,j}}(x)-|x_j|
	=\sqrt{x_j^2+\varepsilon^2|x|^2}-|x_j|
	\le\varepsilon|x|.
	\]
	Thus, for $u\in\mathbb S^{n-1}$, we have
	\[
	\begin{aligned}
		0&\le h_M(u)-h_C(u)
		=h_{A_{k,T}(\lambda)}(D^{-1}u)-\sum_{j=1}^n|u_j|\\
		&=\sum_{j=1}^k
		\bigl(h_{E_{\varepsilon,j}}(D^{-1}u)-|(D^{-1}u)_j|\bigr)
		+\lambda T\sum_{j=k+1}^n
		\bigl(h_{E_{\varepsilon,j}}(D^{-1}u)-|(D^{-1}u)_j|\bigr)\\
		&\le \varepsilon|D^{-1}u|\bigl[k+\lambda T(n-k)\bigr]
		\le \frac{k}{T^2}+\frac{\lambda(n-k)}{T}
		\longrightarrow0,\quad \text{as }T\to+\infty.
	\end{aligned}
	\]
	Hence $M\to C$ in the Hausdorff metric as $T\to\infty$.
	
	For $u\in\mathbb S^{n-1}$, put $r=|D^{-1}u|$ and
	$v=\frac{D^{-1}u}{r}$, so that $Dv=\frac{u}{r}$. Let $P_u$ be the orthogonal
	projection onto $u^\perp$. For $x\in v^\perp$,
	\[
	D(x+\mathbb Rv)=P_uDx+\mathbb Ru,\qquad
	\ell_{DM}(P_uDx,u)=r^{-1}\ell_M(x,v).
	\]
	In orthonormal bases ending in $v$ and $u$, the matrix of $D$ is
	block triangular, with diagonal blocks $P_uD|_{v^\perp}$ and $r^{-1}$.
	Therefore
	\[
	\operatorname{Jac}_{n-1}(P_uD|_{v^\perp})=r|\det D|,
	\qquad
	Z_q(DM,u)=|\det D|\,r^{1-q}Z_q(M,v).
	\]
	
	Write $u^{(k)}=(u_1,\ldots,u_k,0,\ldots,0)$. Then as $T\to+\infty$, we have
	\[
	D^{-1}u\longrightarrow u^{(k)},
	\qquad\text{and}\qquad
	\frac{D^{-1}u}{|D^{-1}u|}
	\longrightarrow \frac{u^{(k)}}{|u^{(k)}|}
	\ \text{if }u^{(k)}\ne0.
	\]
	
	By Lemma~\ref{lem:continuity}, $Z_q(M,\cdot)$ is uniformly bounded
	for large $T$.  By $r\le1$, $1-q>0$  and the dominated convergence theorem applied to
	\eqref{eq:chord-directional}, we have
	\begin{equation}\label{negc:eq:box-limit}
		\frac{I_q(A_{k,T}(\lambda))}{(\lambda T)^{n-k}}
		\longrightarrow c_k \qquad \text{as } T\to+\infty,
	\end{equation}
	where
	$
	c_k=\frac1{n\omega_n}
	\int_{\{u\in\mathbb S^{n-1}:u^{(k)}\ne0\}}
	|u^{(k)}|^{1-q}
	Z_q\!\bigl(C,\frac{u^{(k)}}{|u^{(k)}|}\bigr)\dd u>0
	$
	is independent of $\lambda$.
	
	Let $\alpha=\frac{1}{n+q-1}$ and $\beta_k=(n-k)\alpha$.
	By \eqref{negc:eq:box-limit} and the $1$-homogeneity of $I_q^\alpha$,
	\[
	\frac{I_q(K_k+L_k)^\alpha-I_q(K_k)^\alpha-I_q(L_k)^\alpha}
	{T^{\beta_k}}
	\longrightarrow
	c_k^\alpha\bigl[2\bigl(\frac{3}{2}\bigr)^{\beta_k}-1-2^{\beta_k}\bigr] \qquad \text{as } T\to+\infty.
	\]
	Since $-1<q<0$ and $n\ge3$, we have
	\[
	\begin{array}{lll}
		k=1: & \beta_1=\frac{n-1}{n+q-1}>1,       & 2\bigl(\frac{3}{2}\bigr)^{\beta_1}<1+2^{\beta_1},\\
		k=2: & 0<\beta_2=\frac{n-2}{n+q-1}<1,     & 2\bigl(\frac{3}{2}\bigr)^{\beta_2}>1+2^{\beta_2}.
	\end{array}
	\]
	Thus both required strict inequalities hold for all sufficiently large $T$.
\end{proof}

\vskip 25pt
\section{\bfseries \texorpdfstring{$0\leq q<1$}{0 $\leq$ q < 1}}
\label{sec:fractional-lifting}

\begin{theorem}
	\label{lift:thm:weighted}
	Let $n\ge2$, $0<q<1$, and $\nu$ be a nonzero finite Borel measure
	on $\mathbb S^{n-1}$. Then for $K,L\in\mathcal K^n$, 
	\begin{equation}\label{lift:eq:weighted-bm}
		J_{q,\nu}(K+L)^{\frac1{n+q-1}}
		\ge J_{q,\nu}(K)^{\frac1{n+q-1}}+J_{q,\nu}(L)^{\frac1{n+q-1}}.
	\end{equation}
	If $n=2$ or $\operatorname{supp}\nu=\mathbb S^{n-1}$, equality holds
	if and only if $K$ and $L$ are homothetic.
\end{theorem}
\begin{proof}
	For $n=2$, the result follows from Theorem~\ref{negp:thm:reverse}.
	For $n\ge3$, let $K_t=(1-t)K+tL$ and $\alpha=\frac{1}{n+q-1}$.
	By Lemma~\ref{lem:homogeneous-bm}, to prove \eqref{lift:eq:weighted-bm}
	it suffices to show that
	\begin{equation}\label{lift:eq:segment-bm}
		J_{q,\nu}(K_t)^\alpha\ge (1-t)J_{q,\nu}(K)^\alpha+tJ_{q,\nu}(L)^\alpha,
		\qquad 0<t<1.
	\end{equation}
	
	\textbf{Case 1}. $\nu$ is a Discrete measures.
	
	\emph{(i) At most two atoms.}
	Let $\nu=a\delta_u+b\delta_v$, where $a,b\ge0$ and $a+b>0$, and choose
	a two-dimensional subspace $E$ containing $u$ and $v$. For
	$A\in\{K,L,K_t\}$ and $z\in E^\perp$, let $A_z=\{x\in E:x+z\in A\}$. Then $\mathrm{int} A_z \neq \emptyset$, for $z\in\operatorname{int}(A|E^\perp)$. Let
	\[
	f_A(z)=aZ_q^{(2)}(A_z,u)+bZ_q^{(2)}(A_z,v)
	\qquad\text{for }z\in\operatorname{int}(A|E^\perp),
	\]
	and $f_A=0$ for $z\notin\operatorname{int}(A|E^\perp)$. By Fubini's theore, we have
	$\int_{E^\perp}f_A(z)\dd z=J_{q,\nu}(A)$.
	
	Let $z_t=(1-t)z_0+tz_1$. If $f_K(z_0)f_L(z_1)>0$, then
	$z_t\in\operatorname{int}(K_t|E^\perp)$ and
	$(1-t)K_{z_0}+tL_{z_1}\subseteq(K_t)_{z_t}$. By this inclusion, the
	monotonicity of $Z_q^{(2)}$ and Theorem~\ref{negp:thm:reverse}, we have
	\[
	f_{K_t}(z_t)^{\frac1{q+1}}\ge
	(1-t)f_K(z_0)^{\frac1{q+1}}+tf_L(z_1)^{\frac1{q+1}}.
	\]
	By the Borell--Brascamp--Lieb inequality \cite{BL76} in $E^\perp$,
	with parameter $\frac{1}{q+1}$ and the mean taken to be zero if either entry
	is zero, we obtain \eqref{lift:eq:segment-bm} with exponent
	$\alpha=\frac{1}{q+1+n-2}$.
	
	\emph{(ii) Finitely many atoms, at least three.}
	Let $\nu=\sum_{i=1}^m c_i\delta_{u_i}$ with $c_i>0$ and $m\ge 3$. Consider the extremal 
	\[
	\min\big\{\sum_{i=1}^m p_i Z_q(K_t,u_i)\colon
	p_i\ge0,\ \sum_{i=1}^m p_i Z_q(A,u_i)=J_{q,\nu}(A)
	\text{ for } A=K,L\big\}.
	\]
	Since
	$0\le p_i\le \frac{J_{q,\nu}(K)}{Z_q(K,u_i)}$, the feasible set contains $(c_i)_{i=1}^m$ and is compact. Let $\bar{p}$ be a minimizer. If at least three coordinates of $\bar{p}$ are positive, 
	there is a nonzero vector $b$ supported on these coordinates with
	$\sum_{i=1}^m b_iZ_q(A,u_i)=0$ for $A=K,L$. Then $\bar{p}\pm\varepsilon b$
	are feasible for sufficiently small $\varepsilon>0$, contradicting
	extremality. Thus the measure $\mu=\sum_{i=1}^m \bar{p}_i\delta_{u_i}$ is supported on
	at most two atoms and is nonzero by the endpoint constraints.
	By minimality, the inequality in (i) and the endpoint constraints,
	we have
	\[
	J_{q,\nu}(K_t)^\alpha\ge J_{q,\mu}(K_t)^\alpha
	\ge (1-t)J_{q,\mu}(K)^\alpha+tJ_{q,\mu}(L)^\alpha
	=(1-t)J_{q,\nu}(K)^\alpha+tJ_{q,\nu}(L)^\alpha.
	\]
	
	\emph{(iii) Countably infinitely many atoms.}
	Write $\nu=\sum_{i=1}^\infty c_i\delta_{u_i}$ and let
	$\nu_m=\sum_{i=1}^m c_i\delta_{u_i}$. By monotone convergence,
	$J_{q,\nu_m}(A)\to J_{q,\nu}(A)$ as $m\to+\infty$ for
	$A\in\{K,L,K_t\}$. Taking limits in \eqref{lift:eq:segment-bm}
	for $\nu_m$ gives the inequality for $\nu$.
	
	\textbf{Case 2}. $\nu$ is a general measures.
	
	We first prove equality \eqref{lift:eq:weighted-bm}.
	For arbitrary $\nu$, choose finite positive atomic measures $\nu_j$
	converging weakly to $\nu$. By continuity of $Z_q(A,\cdot)$, we have
	$J_{q,\nu_j}(A)\to J_{q,\nu}(A)$ as $j\to+\infty$ for
	$A\in\{K,L,K_t\}$. By Case~1 and passage to the limit, we obtain
	\eqref{lift:eq:segment-bm}.
	
	Assume $\mathrm{supp}\,\nu=\mathbb{S}^{n-1}$. 	In the following, we  characterize the equality condition. 
	
	If  $K$ and $L$ are homothetic, the equality in \eqref{lift:eq:weighted-bm} holds by the translation
	invariance and homogeneity of $=J_{q,\nu}$.
	
	Conversely, suppose equality holds in \eqref{lift:eq:weighted-bm}.
	Let $a=J_{q,\nu}(K)^\alpha$ and $b=J_{q,\nu}(L)^\alpha$. Then
	\[
	J_{q,\nu}\bigl(\frac{a}{a+b}\frac{K}{a}+\frac{b}{a+b}\frac{L}{b}\bigr)^\alpha
	=\frac{J_{q,\nu}(K+L)^\alpha}{a+b}=1.
	\]
	By Lemma~\ref{lem:homogeneous-bm},
	we have $J_{q,\nu}(K_t)=1$ for every $t\in[0,1]$.
	For $u\in\mathbb S^{n-1}$ and $\varepsilon>0$, applying
	\eqref{lift:eq:segment-bm} to $\nu+\varepsilon\delta_u$,
	we have
	\[
	(1+\varepsilon Z_q(K_t,u))^\alpha
	\ge (1-t)(1+\varepsilon Z_q(K,u))^\alpha
	+t(1+\varepsilon Z_q(L,u))^\alpha.
	\]
	Subtracting $1$, dividing by $\alpha\varepsilon$ and letting
	$\varepsilon\to0^+$, it follows that $\Delta_t\ge0$, where $\Delta_t(u)=Z_q(K_t,u)-(1-t)Z_q(K,u)-tZ_q(L,u)$.
	Since $\Delta_t\in C\bigl(\mathbb S^{n-1}\times[0,\infty)\bigr)$,  $\mathrm{supp}\,\nu=\mathbb{S}^{n-1}$ and
	\[
	\int_{\mathbb S^{n-1}}\Delta_t(u)\dd\nu(u)=
	=J_{q,\nu}(K_t)-(1-t)J_{q,\nu}(K)-tJ_{q,\nu}(L)=1-(1-t)-t=0,
	\]
	it follows that
	\[
	Z_q(K_t,u)=(1-t)Z_q(K,u)+tZ_q(L,u)
	\qquad\text{for }u\in\mathbb S^{n-1},\ t\in[0,1].
	\]
	Since $\int_{\mathbb S^{n-1}}(Z_q(K,u)-Z_q(L,u))\dd\nu(u)=0$,
	there exists
	$u_0\in\mathbb{S}^{n-1}$ such that $Z_q(K,u_0)=Z_q(L,u_0)$. Hence $Z_q(K_t,u_0)$ is
	constant on $[0,1]$. In the following, we show that $|K_t|$ is constant on $[0,1]$.
	
	For arbitrary $t,s\in[0,1]$, let $f_0=\ell_{K_t}(\cdot,u_0)$,
	$f_1=\ell_{K_s}(\cdot,u_0)$ and $h=\ell_{K_{\frac{t+s}{2}}}(\cdot,u_0)$
	on the interiors of the respective projections onto $u_0^\perp$,
	and extend them by zero. By the constancy of $Z_q(K_t,u_0)$ on $[0,1]$,
	we have $\int f_0^q=\int f_1^q=\int h^q>0$.
	By $K_{\frac{t+s}{2}}=\frac{K_t+K_s}{2}$ and the chord inclusion, we have
	$h\bigl(\frac{x+y}{2}\bigr)\ge\frac{f_0(x)+f_1(y)}{2}$ if $f_0(x)f_1(y)>0$.
	Let $T$ be the Brenier map with $T_\#(f_0^q\dd x)=f_1^q\dd x$,
	and let $S=\frac{\mathrm{Id}+T}{2}$. On $U=\{x:f_0(x)>0\}$, let
	$a(x)=\frac{f_1(Tx)}{f_0(x)}$, which is positive almost everywhere.
	By convexity of the Brenier potential and the Monge--Amp\`ere identity
	\cite[Remark~4.5]{McCann}, $DT$ is positive semidefinite and
	$\det DT(x)=a(x)^{-q}$ almost everywhere on $U$, with derivatives
	in the Alexandrov sense. If $\lambda_j(x)$ are the eigenvalues
	of $DT(x)$, then
	\[
	\det DS(x)=\prod_{j=1}^{n-1}\frac{1+\lambda_j(x)}2
	\ge\sqrt{\det DT(x)}=a(x)^{-\frac{q}{2}}.
	\]
	Since $S$ is a gradient of a convex function, by the change of
	variables \cite[Corollary~4.3]{McCann}, chord inclusion and the
	arithmetic--geometric mean inequality, we have
	\[
	\begin{aligned}
		\int h^q
		&\ge\int_U h(Sx)^q\det DS(x)\dd x
		\ge\int_U f_0(x)^q\bigl(\frac{1+a(x)}2\bigr)^q\det DS(x)\dd x\\
		&\ge\int_U f_0(x)^q\bigl(\frac{1+a(x)}{2\sqrt{a(x)}}\bigr)^q\dd x
		\ge\int f_0^q=\int h^q.
	\end{aligned}
	\]
	Equality in the last inequality gives $a=1$ almost everywhere for
	$f_0^q\dd x$. By the pushforward identity, we have
	\[
	|K_s|=\int f_1=\int_U f_1(Tx)^{1-q}f_0(x)^q\dd x=\int f_0=|K_t|.
	\]
	Thus $|K_t|$ is constant. By the equality case of the
	Brunn--Minkowski inequality for volume, $K$ and $L$
	are translates. Rescaling gives homothety.
\end{proof}

For $0<q<1$, taking $\dd\nu(u)=\frac{\dd u}{n\omega_n}$, we obtain
$J_{q,\nu}=I_q$. Thus Theorem~\ref{lift:thm:weighted} proves
Theorem~\ref{thm:chord-main} in this range.

For $q=0$, by \eqref{intro:eq:classical-values} and the Brunn--Minkowski
inequality for surface area \cite[Chapter~7]{Schneider}, we have
\[
I_0(K+L)^{\frac1{n-1}}
\ge I_0(K)^{\frac1{n-1}}+I_0(L)^{\frac1{n-1}}.
\]
For $n\ge3$, the equality case of the surface area inequality gives
equality if and only if $K$ and $L$ are homothetic. For $n=2$, since
$S(K)=\int_{\mathbb S^1}h_K(u)\dd u$ is Minkowski linear,  equality
holds for all $K,L$.

\vskip 25pt
\section{\bfseries \texorpdfstring{$1\le q\le n+1$}{1 <= q <= n+1}}
\label{sec:covariogram}

\begin{theorem}\label{thm:energy}
	For $K,L\in\mathcal{K}^n$ and $0\le a<n$,
	\begin{equation}\label{eq:energy-bm}
		E_a(K+L)^{\frac1{2n-a}}
		\ge E_a(K)^{\frac1{2n-a}}+E_a(L)^{\frac1{2n-a}},
	\end{equation}
	with equality if and only if $K$ and $L$ are homothetic.
\end{theorem}

\begin{proof}
	Let $K_t=(1-t)K+tL$, for $0<t<1$, $p=2n-a$,
	$g_K(x)=|K\cap(K+x)|$ and
	\[
	C_K=\{(x,y)\in\R^{2n}:x\in K-K,\ |y|\le g_K(x)^{\frac1n}\}.
	\]
	For $x\in K-K$ and $z\in L-L$, we have
	\[
	(1-t)(K\cap(K+x))+t(L\cap(L+z))
	\subseteq K_t\cap\bigl(K_t+(1-t)x+tz\bigr).
	\]
	By the Brunn--Minkowski inequality for volume, we have
	\[
	g_{K_t}\bigl((1-t)x+tz\bigr)^{\frac1n}
	\ge (1-t)g_K(x)^{\frac1n}+t g_L(z)^{\frac1n}.
	\]
	Taking $L=K$ in this inequality shows that $g_K^{\frac1n}$ is concave
	on $K-K$. By the preceding inequality and the triangle inequality,
	$C_K$ is convex and $(1-t)C_K+tC_L\subseteq C_{K_t}$.
	Since $g_K$ is continuous and even and $g_K(0)=|K|>0$, it follows
	that $C_K$ is compact, origin-symmetric and has nonempty interior.
	
	Let
	$\dd\mu_a(x,y)=|x|^{-a}\dd x\dd y$. By Tonelli's theorem, we have
	\begin{align*}
		\mu_a(C_K)
		&=\int_{C_K}|x|^{-a}\dd x\dd y
		=\int_{K-K}|x|^{-a}
		\bigl(\int_{|y|\le g_K(x)^{\frac{1}{n}}}\dd y\bigr)\dd x\\
		&=\omega_n\int_{K-K}|x|^{-a}g_K(x)\dd x
		=\omega_n\int_{K-K}\int_{K\cap(K+x)}|x|^{-a}\dd z\dd x\\
		&=\omega_n\int_K\int_K|z-w|^{-a}\dd w\dd z
		=\omega_n E_a(K).
	\end{align*}
	By $(1-t)C_K+tC_L\subseteq C_{K_t}$, we have
	$\mu_a(C_{K_t})\ge\mu_a\bigl((1-t)C_K+tC_L\bigr)$. To prove inequality \ref{eq:energy-bm}, it suffices to prove that
	\begin{equation}\label{eq:energy-lift-weighted}
		\mu_a\bigl((1-t)C_K+tC_L\bigr)^{\frac1p}
		\ge (1-t)\mu_a(C_K)^{\frac1p}+t\mu_a(C_L)^{\frac1p}.
	\end{equation}
	
	For $a=0$, \eqref{eq:energy-lift-weighted} is the Brunn--Minkowski
	inequality for volume in $\R^{2n}$, since $\mu_0$ is volume and $p=2n$.
	
	For $a>0$, applying Sadovsky and Zhang's inequality \cite[Theorem~2]{SZ}
	for Minkowski combinations in dimension $2n$ with index $p=2n-a$, it follows that
	\begin{align*}
		&\ \ \ \   \mu_a\bigl((1-t)C_K+tC_L\bigr)^{\frac1p}=\lim_{\varepsilon\to0^+}
		\bigl[\int_{(1-t)C_K+tC_L}
		\bigl(|x|^2+\varepsilon^2|y|^2\bigr)^{-\frac{a}{2}}
		\dd x\dd y\bigr]^{\frac1p}\\
		&=\lim_{\varepsilon\to0^+}\varepsilon^{-\frac{n}{p}}
		\bigl[\int_{\{(x,\varepsilon y):(x,y)\in(1-t)C_K+tC_L\}}
		|z|^{-a}\dd z\bigr]^{\frac1p}\\
		&\ge\lim_{\varepsilon\to0^+}\varepsilon^{-\frac{n}{p}}
		\bigl\{(1-t)\bigl[\int_{\{(x,\varepsilon y):(x,y)\in C_K\}}
		|z|^{-a}\dd z\bigr]^{\frac1p}
		+t\bigl[\int_{\{(x,\varepsilon y):(x,y)\in C_L\}}
		|z|^{-a}\dd z\bigr]^{\frac1p}\bigr\}\\
		&=\lim_{\varepsilon\to0^+}\bigl\{(1-t)\bigl[\int_{C_K}
		\bigl(|x|^2+\varepsilon^2|y|^2\bigr)^{-\frac{a}{2}}
		\dd x\dd y\bigr]^{\frac1p}+t\bigl[\int_{C_L}
		\bigl(|x|^2+\varepsilon^2|y|^2\bigr)^{-\frac{a}{2}}
		\dd x\dd y\bigr]^{\frac1p}\bigr\}\\
		&=(1-t)\mu_a(C_K)^{\frac1p}+t\mu_a(C_L)^{\frac1p}.
	\end{align*}
	Here $|x|^{-a}$ is integrable on bounded subsets of $\R^{2n}$ since
	$a<n$. This proves \eqref{eq:energy-lift-weighted}.
	
	If $K,L$ are homothetic, equality in \eqref{eq:energy-bm} follows from
	translation invariance and homogeneity of $E_a$. Conversely, suppose
	equality holds in \eqref{eq:energy-bm}. By Lemma~\ref{lem:homogeneous-bm},
	equality also holds along $K_t$. By $(1-t)C_K+tC_L\subseteq C_{K_t}$,
	\eqref{eq:energy-lift-weighted}, and $\mu_a(C_K)=\omega_n E_a(K)$, we have
	\begin{align*}
		\mu_a(C_{K_t})^{\frac1p}
		&\ge \mu_a\bigl((1-t)C_K+tC_L\bigr)^{\frac1p}
		\ge (1-t)\mu_a(C_K)^{\frac1p}+t\mu_a(C_L)^{\frac1p}=\mu_a(C_{K_t})^{\frac1p}.
	\end{align*}
	Since
	$\mu_a$ has positive density almost everywhere, we have
	$C_{K_t}=(1-t)C_K+tC_L$.
	
	From that
	$C_{K_t}\cap\bigl(\{0\}\times\R^n\bigr)=\{0\}\times |K_t|^{\frac1n}B^n$, $g_K\le|K|$, $g_L\le|L|$,
	and the Brunn--Minkowski inequality for volume, it follows that
	\begin{align*}
		|K_t|^{\frac1n}
		&\le\max_{(x,y)\in(1-t)C_K+tC_L}|y|\le(1-t)|K|^{\frac1n}+t|L|^{\frac1n}
		\le |K_t|^{\frac1n}.
	\end{align*}
	By the equality condition of the Brunn--Minkowski inequality for volume, it follows that $K,L$ are homothetic.
\end{proof}

At $q=1$, the inequality is the classical Brunn--Minkowski inequality
for volume.

For $1<q\le n+1$, taking $a=n+1-q$, by \eqref{intro:eq:chord-energy} and Theorem~\ref{thm:energy}, we obtain
the Brunn--Minkowski inequality for $I_q$, with equality if and only
if $K$ and $L$ are homothetic. Together with $I_1=|K|$, this proves
Theorem~\ref{thm:chord-main} for $1\le q\le n+1$.

\vskip 25pt
\section{\bfseries \texorpdfstring{$n+1<q\le n+3+\eta_n$}{n+1 < q <= n+3+eta_n}}
\label{sec:cnd-full-interval}

For a finite Borel measure $\nu$ on $\mathbb S^{n-1}$, $p>0$ and
$q>n+1$, define
\begin{align*}
 k_{p,\nu}(z)&=\int_{\mathbb S^{n-1}}|\theta\cdot z|^p\dd\nu(\theta),\\
 E_{q,\nu}(K)&=\int_{K^2}k_{q-n-1,\nu}(x-y)\dd x\dd y.
\end{align*}

We prove the Brunn--Minkowski inequality for $I_q$ by establishing it for $E_{q,\nu}$.

\begin{theorem}\label{aw:thm:angular-window}
There exists a constant $0<\eta_n\le\frac14$,
depending only on $n$, such that if $n+1<q\le n+3+\eta_n$ and
$\nu$ is a nonzero finite Borel measure on $\mathbb S^{n-1}$,
then for $K,L\in\mathcal K^n$,
\[
 E_{q,\nu}(K+L)^{\frac1{n+q-1}}
 \ge E_{q,\nu}(K)^{\frac1{n+q-1}}+E_{q,\nu}(L)^{\frac1{n+q-1}},
\]
with equality if and only if $K$ and $L$ are homothetic.
\end{theorem}

\subsection{\texorpdfstring{$n+1<q\le n+3$}{n+1 < q <= n+3}}

\begin{lemma}\label{lem:one-dimensional-square}
	If $\dd\mu(x)=e^{-\varphi(x)}\dd x$ is a probability measure on $\R$
	with $\varphi\in C^2(\R)$ and $0<c\le\varphi''\le C$, then for
	$0<r\le2$ and $f\in C^1(\R)$ with bounded derivative and $\int_\R f\dd\mu=0$,
	\begin{equation}
		\label{eq7.1} \frac12\int_{\R^2}|x-y|^r\bigl(f(x)+f(y)\bigr)^2\dd\mu(x)\dd\mu(y)
		\le \int_{\R^2}|x-y|^r\frac{f'(x)^2}{\varphi''(x)}\dd\mu(x)\dd\mu(y).
	\end{equation}
\end{lemma}
\begin{proof}
	For Borel functions $a,b:\R\to\R$ of  with bounded derivative , define
	\[
	\mathcal C(a,b):=\frac12\int_{\R^2}|x-y|^r
	\bigl(a(x)+a(y)\bigr)\bigl(b(x)+b(y)\bigr)\dd\mu(x)\dd\mu(y).
	\]
	The integral is finite since $\varphi''\ge c>0$ implies that $\mu$
	has Gaussian tails.
	
	For $s\in\R$, let
	$g_s(x)=\mathbf1_{\{x>s\}}-\mu((s,\infty))$.
	
	\textbf{Step 1.} Write $K(s,t)=\mathcal C(g_s,g_t)$, for $s,t\in\R$. In the following,
	We prove that $K(s,t)\ge0$. By symmetry, it suffices to consider $s\le t$.
	
	Let $\alpha=\mu((-\infty,s])$ and $\beta=\mu((t,\infty))$.
	For every Borel set $\omega\subseteq\R$, define
	\[
	\mu_-(\omega):=\frac{\mu\bigl(\omega\cap(-\infty,s]\bigr)}{\alpha},
	\qquad
	\mu_+(\omega):=\frac{\mu\bigl(\omega\cap(t,\infty)\bigr)}{\beta}.
	\]
	Then $ g_s=\alpha-\mathbf1_{(-\infty,s]},$ and $ g_t=\mathbf1_{(t,\infty)}-\beta.$
	
	For a finite signed measure $\sigma$ on $\R$ with finite absolute
	$r$-th moment, let
	\[
	\mathcal Q(\sigma)=\int_{\R^2}|x-y|^r\dd\sigma(x)\dd\sigma(y).
	\]
	By $(-\infty,s]\cap(t,\infty)=\varnothing$,
	we have
	\begin{align*}
		\frac{2K(s,t)}{\alpha\beta}
		&=\int_{\R^2}|x-y|^r
		\bigl(2-\frac{\mathbf1_{\{x\le s\}}+\mathbf1_{\{y\le s\}}}{\alpha}\bigr)
		\bigl(\frac{\mathbf1_{\{x>t\}}+\mathbf1_{\{y>t\}}}{\beta}-2\bigr)
		\dd\mu(x)\dd\mu(y)\\
		&=-4\mathcal Q(\mu)
		+4\int_{\R^2}|x-y|^r\dd(\mu_-+\mu_+)(x)\dd\mu(y)-2\int_{\R^2}|x-y|^r\dd\mu_-(x)\dd\mu_+(y)\\
		&=\mathcal Q(\mu_-)+\mathcal Q(\mu_+)
		-\mathcal Q(\mu_-+\mu_+-2\mu).
	\end{align*}
	
	Assume $\sigma(\R)=0$. For $r=2$, we have
	$\mathcal Q(\sigma)=-2\bigl(\int_\R x\dd\sigma(x)\bigr)^2\le0$.
	For $0<r<2$, there exists $c_r>0$ such that
	\[
	|z|^r=c_r\int_0^\infty\bigl(1-\cos(\xi z)\bigr)
	\frac{\dd\xi}{\xi^{1+r}}.
	\]
	By the moment assumption of $\sigma$, $\sigma(\R)=0$ and Fubini's theorem, we have
	\[
	\mathcal Q(\sigma)=-c_r\int_0^\infty
	\biggl|\int_\R e^{i\xi x}\dd\sigma(x)\biggr|^2
	\frac{\dd\xi}{\xi^{1+r}}\le0.
	\]
	Since $\mathcal Q(\mu_\pm)\ge0$ and
	$(\mu_-+\mu_+-2\mu)(\R)=1+1-2=0$,
	it follows that $K(s,t)\ge0$.

	\vskip 3pt
	\textbf{Step 2.} Prove inequality \eqref{eq7.1}.
	
	For $b\in C^1(\R)$ with bounded derivative, by integrating
	$b(x)-b(y)=\int_y^x b'(s)\dd s$ with respect to $\mu(y)$, we have
	\[
	b(x)-\int_\R b\dd\mu=\int_\R b'(s)g_s(x)\dd s.
	\]
	Taking $b=f$ and using $\int_\R f\dd\mu=0$, we have
	$f(x)=\int_\R f'(s)g_s(x)\dd s$. Substituting this identity into
	$\mathcal C(f,f)$, we obtain
	\[
	\mathcal C(f,f)
	=\frac12\int_{\R^2}|x-y|^r
	\int_{\R^2}f'(s)f'(t)
	\bigl(g_s(x)+g_s(y)\bigr)\bigl(g_t(x)+g_t(y)\bigr)
	\dd s\dd t\dd\mu(x)\dd\mu(y).
	\]
	Since $\mathbf1_{\{x>s\}}$ takes values in $\{0,1\}$, by the definition
	of $g_s$ and Tonelli's theorem, we have
	\[
	\int_\R|g_s(x)|\dd s
	=\int_\R\int_\R
	\bigl|\mathbf1_{\{x>s\}}-\mathbf1_{\{y>s\}}\bigr|
	\dd s\dd\mu(y)
	=\int_\R|x-y|\dd\mu(y).
	\]
	Thus,
	\begin{align*}
		&\int_{\R^2}|x-y|^r\int_{\R^2}
		|f'(s)f'(t)|\,|g_s(x)+g_s(y)|\,|g_t(x)+g_t(y)|
		\dd s\dd t\dd\mu(x)\dd\mu(y)\\
		\le& \|f'\|_\infty^2
		\int_{\R^2}|x-y|^r
		\bigl(|x|+|y|+2\int_\R |z|\dd\mu(z)\bigr)^2
		\dd\mu(x)\dd\mu(y)<\infty.
	\end{align*}
	The last integral is finite by the Gaussian tails of $\mu$. This
	absolute integrability allows us to apply Fubini's theorem to integrate
	in $x,y$ first. By the definition of $K(s,t)$, we have
	\begin{align*}
		\mathcal C(f,f)
		&=\int_{\R^2} f'(s)f'(t)
		\bigl[\frac12\int_{\R^2}|x-y|^r
		\bigl(g_s(x)+g_s(y)\bigr)\bigl(g_t(x)+g_t(y)\bigr)
		\dd\mu(x)\dd\mu(y)\bigr]\dd s\dd t\\
		&=\int_{\R^2}K(s,t)f'(s)f'(t)\dd s\dd t.
	\end{align*}
	
	By the decay of $e^{-\varphi}$ and the integral representation with
	$b=\varphi'$, we have
	\begin{align*}
		\int_\R\varphi'\dd\mu
		&=-\int_\R(e^{-\varphi(x)})'\dd x
		=-\bigl[e^{-\varphi(x)}\bigr]_{-\infty}^{+\infty}=0,\\
		\varphi'(x)
		&=\varphi'(x)-\int_\R\varphi'\dd\mu
		=\int_\R\varphi''(t)g_t(x)\dd t.
	\end{align*}
	By $|g_s|\le1$, $\varphi''\le C$, the formula for
	$\int_\R|g_t(x)|\dd t$ and the Gaussian tails of $\mu$, we have
	\begin{align*}
		&\frac12\int_{\R^2}|x-y|^r\int_\R
		\bigl|\varphi''(t)(g_s(x)+g_s(y))(g_t(x)+g_t(y))\bigr|
		\dd t\dd\mu(x)\dd\mu(y)\\
		&\qquad\le C\int_{\R^2}|x-y|^r
		\biggl(|x|+|y|+2\int_\R|z|\dd\mu(z)\biggr)\dd\mu(x)\dd\mu(y)<\infty.
	\end{align*}
	By Fubini's theorem, the substitution $z=y-x$ and integration by parts
	with $\dd g_s=\delta_s$, we have
	\begin{align*}
		&\ \ \ \ \int_\R K(s,t)\varphi''(t)\dd t=\mathcal C(g_s,\varphi')\\
		&=\int_{\R^2}|x-y|^r g_s(x)
		\bigl(\varphi'(x)+\varphi'(y)\bigr)\dd\mu(x)\dd\mu(y)\\
		&=-\int_\R |z|^r\int_\R g_s(x)\,\partial_x
		\bigl(e^{-\varphi(x)-\varphi(x+z)}\bigr)\dd x\dd z\\
		&=\int_\R |z|^r e^{-\varphi(s)-\varphi(s+z)}\dd z
		=e^{-\varphi(s)}\int_\R|s-y|^r\dd\mu(y).
	\end{align*}
	
	Observe that
	\[
	2f'(s)f'(t)
	\le\frac{f'(s)^2}{\varphi''(s)}\varphi''(t)
	+\frac{f'(t)^2}{\varphi''(t)}\varphi''(s).
	\]
	Multiplying by $K(s,t)\ge0$, integrating and using symmetry of $K$,
	we obtain
	\begin{align*}
		\mathcal C(f,f)
		\le\int_\R\frac{f'(s)^2}{\varphi''(s)}
		\biggl(\int_\R K(s,t)\varphi''(t)\dd t\biggr)\dd s=\int_{\R^2}|s-y|^r\frac{f'(s)^2}{\varphi''(s)}
		\dd\mu(s)\dd\mu(y).
	\end{align*}
	This proves \eqref{eq7.1}.
\end{proof}

\begin{proof}[Proof of Theorem~\ref{aw:thm:angular-window} for $n+1<q\le n+3$]
	Write $r=q-n-1\in(0,2]$.
	\vskip 3pt
	\textbf{Step 1.} Prove a differential inequality for smooth densities.
	
	For an open interval $I$, let $\Phi\in C^2(I\times\R^n)$ be jointly convex,
	with $D_x^2\Phi\ge c\,\mathrm{Id}$ and
	$D^2\Phi\le M\,\mathrm{Id}$ for $0<c\le M$. Write $z=\log Z$ and
	$g=\log H$, where
	\[
	Z(t)=\int_{\R^n} e^{-\Phi(t,x)}\dd x,\qquad
	H(t)=\int_{\R^{2n}} k_{r,\nu}(x-y)e^{-\Phi(t,x)-\Phi(t,y)}\dd x\dd y .
	\]
	On compact subintervals of $I$,
	\[
	\Phi(t,x)\ge\frac c4|x|^2-C,\qquad
	|\nabla_{t,x}\Phi(t,x)|\le C(1+|x|).
	\]
	These estimates and the Hessian bound justify differentiation under the integrals.
	
	For $\theta\in\mathbb S^{n-1}$, define the marginal density and its energy by
	\begin{align*}
		e^{-V(t,s)}&=\int_{\theta^\perp}e^{-\Phi(t,s\theta+y)}\dd y,\\
		H_\theta(t)&=\int_{\R^2}|s-y|^r e^{-V(t,s)-V(t,y)}\dd s\dd y.
	\end{align*}
	By Fubini's theorem and the definition of $k_{r,\nu}$,
	\[
	Z(t)=\int_\R e^{-V(t,s)}\dd s,\qquad
	H(t)=\int_{\mathbb S^{n-1}}H_\theta(t)\dd\nu(\theta).
	\]
	
	By Pr\'ekopa's theorem \cite{BL76}, $V$ is jointly convex. Applying it to
	$\Phi(t,s\theta+y)-\frac c2s^2$ at fixed $t$ yields $V_{ss}\ge c$.
	For a unit direction $a$ in the $(t,s)$-plane, by differentiating
	the marginal integral twice and applying Cauchy--Schwarz, we have
	\[
	\partial_{aa}V(t,s)\le
	\frac{\displaystyle\int_{\theta^\perp}
		\partial_{aa}\bigl[\Phi(t,s\theta+y)\bigr]e^{-\Phi(t,s\theta+y)}\dd y}
	{\displaystyle\int_{\theta^\perp}e^{-\Phi(t,s\theta+y)}\dd y}
	\le M.
	\]
	
	For fixed $t$, let $\dd\mu(s)=e^{-V(t,s)-z(t)}\dd s$ and
	$f(s)=-V_t(t,s)-z'$. By differentiating $Z$, we have
	$z'=-\int_\R V_t(t,s)\dd\mu(s)$, so $f$ has mean zero.
	By the Hessian bounds, we have
	$|f'|=|V_{ts}|\le M$ and
	$(f')^2\le V_{tt}V_{ss}$. Differentiating $H_\theta$ twice and
	applying Lemma~\ref{lem:one-dimensional-square} with $\varphi=V(t,\cdot)+z(t)$, we obtain
	\begin{align*}
		\frac{H_\theta''-4z'H_\theta'+4(z')^2H_\theta}{Z^2}
		&=\int_{\R^2}|s-y|^r\bigl(f(s)+f(y)\bigr)^2\dd\mu(s)\dd\mu(y)\\
		&\quad-2\int_{\R^2}|s-y|^r V_{tt}(t,s)\dd\mu(s)\dd\mu(y)\le0.
	\end{align*}
	Combining this with $H=\int_{\mathbb S^{n-1}}H_\theta\dd\nu(\theta)>0$ and $g=\log H$, it follows that
	\begin{align}
		\label{eq:cnd-functional-deficit}g''+(g'-2z')^2
		&=\frac{H''}{H}-4z'\frac{H'}{H}+4(z')^2
		=\frac{H''-4z'H'+4(z')^2H}{H}\\
		&=\frac{Z^2}{H}\int_{\mathbb S^{n-1}}
		\frac{H_\theta''-4z'H_\theta'+4(z')^2H_\theta}{Z^2}\dd\nu(\theta)
		\le0.\notag
	\end{align}
	
	\vskip 3pt
	\textbf{Step 2.} Prove the desired inequality.
	
	Let $K_t=(1-t)K+tL$ for $t\in[0,1]$. By Lemma~\ref{lem:homogeneous-bm}
	applied to $E_{q,\nu}^{\frac1{n+q-1}}$, it suffices to prove
	$E_{q,\nu}(K_t)\ge1$ when $E_{q,\nu}(K)=E_{q,\nu}(L)=1$.
	Choose $R>0$ with $K,L\subseteq RB^n$. The convex hull
	$D=\operatorname{conv}((\{0\}\times K)\cup(\{1\}\times L))$
	has section $K_t$ at height $t$. For a smooth probability mollifier
	$\rho_\delta$ supported in $\delta B^{n+1}$, define
	\[
	\Phi_\delta(t,x)=\frac{(d(\cdot,D)^2*\rho_\delta)(t,x)}{2\delta}
	+\frac{\delta}{2}|x|^2,\qquad 0<\delta\le1 .
	\]
	The squared distance to $D$ is convex with $2$-Lipschitz gradient.
	Thus $\Phi_\delta$ is convex and
	\[
	D_x^2\Phi_\delta\ge\delta\,\mathrm{Id},\qquad
	D^2\Phi_\delta\le\left(\frac1\delta+\delta\right)\mathrm{Id}.
	\]
	For $\zeta=(t,x)$, by the $1$-Lipschitz bound for $d(\cdot,D)$ and
	the support of $\rho_\delta$, we have
	\[
	(d(\zeta,D)-\delta)_+^2
	\le (d(\cdot,D)^2*\rho_\delta)(\zeta)
	\le (d(\zeta,D)+\delta)^2.
	\]
	After division by $2\delta$, the convolution tends to zero on $D$ and to infinity off $D$. Hence $e^{-\Phi_\delta(t,x)}\to\mathbf1_{K_t}(x)$ for $t\in[0,1]$.
	Since $d((t,x),D)\ge(|x|-R)_+$, the same estimate bounds these densities by
	$\exp\bigl(-\frac12(|x|-R-1)_+^2\bigr)$. By this bound,
	$k_{r,\nu}(x-y)\le\nu(\mathbb S^{n-1})|x-y|^r$ and dominated convergence,
	we have
	\[
	g_\delta=\log H_\delta\longrightarrow g=\log E_{q,\nu}(K_t),
	\qquad z_\delta=\log Z_\delta\longrightarrow z=\log|K_t| .
	\]
	By \eqref{eq:cnd-functional-deficit}, $g_\delta$ and hence $g$ are concave.
	Since $g(0)=g(1)=0$, we have $E_{q,\nu}(K_t)=e^{g(t)}\ge1$ on $[0,1]$.
	
	\vskip 3pt
	\textbf{Step 3.} Characterize equality.
	
	Suppose that $E_{q,\nu}(K_{t_0})=1$ for some $0<t_0<1$.
	Since $g(0)=g(t_0)=g(1)=0$ and $g$ is concave, $g\equiv0$.
	For $D_\delta=g_\delta-2z_\delta$ and $0<a<b<1$, by Cauchy--Schwarz,
	\eqref{eq:cnd-functional-deficit} and concavity of $g_\delta$, we have
	\begin{align*}
		\frac{|D_\delta(b)-D_\delta(a)|^2}{b-a}
		&\le\int_a^b(D_\delta')^2\dd t\le g_\delta'(a)-g_\delta'(b)\\
		&\le\frac{g_\delta(a)-g_\delta(0)}a
		-\frac{g_\delta(1)-g_\delta(b)}{1-b}\longrightarrow0.
	\end{align*}
	Since $D_\delta\to-2z$, we have $z(a)=z(b)$. By continuity, $|K_t|$ is
	constant on $[0,1]$. The equality condition in the Brunn--Minkowski
	inequality for volume shows that the normalized bodies are translates,
	so the original bodies are homothetic. The converse follows by homogeneity
	and translation invariance.
\end{proof}

\subsection{\texorpdfstring{$n+3<q\le n+3+\eta_n$}{n+3 < q <= n+3+eta_n}}
\label{sec:supercritical-restored}
\label{aw:sec:angular-window}

\begin{proof}[Proof of Theorem~\ref{aw:thm:angular-window}]
Write $r=q-n-1$. Then $2<r\le\frac94$. Let $K$ be smooth with positive principal
curvatures. In John's normalization
$K=b+T\widehat K$, $B^n\subseteq\widehat K\subseteq nB^n$,
\[
 E_{q,\nu}(K)=|\det T|^2 E_{q,\nu_*}(\widehat K),\qquad
 \nu_*=\left(\theta\mapsto\frac{T^\top\theta}{|T^\top\theta|}\right)_*
                       (|T^\top\theta|^r\dd\nu).
\]
Write $K$ for $\widehat K$ and $\nu$ for the normalization of $\nu_*$ to
a probability measure. Both the affine map and $\nu$ remain fixed during the
variation and for $2\le p\le r$.

\vskip 3pt
\textbf{Step 1.} Compute the second variation.

For $f\in C^\infty(\mathbb S^{n-1})$ and small $|t|$, let $h_{K_t}=h_K+tf$.
Write
$V(t)=|K_t|$, $k_p=k_{p,\nu}$ and $E_p(t)=E_{n+1+p,\nu}(K_t)$.
To estimate $(E_p^{\frac1{2n+p}})''(0)$, define
\[
 Q_p=E_p''-\frac{2n+p-1}{2n+p}\frac{(E_p')^2}{E_p}
    =(2n+p)E_p^{1-\frac1{2n+p}}\bigl(E_p^{\frac1{2n+p}}\bigr)''.
\]
Identify functions on the sphere with their pullbacks to $\partial K$ by the outer unit normal $N$.
By degree-one homogeneity of $E_p^{\frac1{2n+p}}$, adding a multiple of
$h_K$ to $f$ preserves $Q_p$. Since $\int_{\partial K}h_K\dd S=nV$,
we may take $\int_{\partial K}f\dd S=0$. Let $u$ be harmonic in $K$ with
$\partial_Nu=f$ on $\partial K$. Subtracting $a\cdot x$ from $u$, with $a=\frac1V\int_K\nabla u\dd x$,
adds the translation velocity $-a\cdot N$ to $f$. Thus $Q_p$ is unchanged and
$v=\nabla u$ has mean zero on $K$.

For $U$ of class $C^2$ near $K$, define
\[
 B(U)=\left.\frac{d^2}{dt^2}\right|_{t=0}\int_{K_t}U(x)\dd x.
\]
With $\mathrm{II}=dN$, $H=\operatorname{tr}\mathrm{II}$ and $\nabla_\partial$
the boundary gradient, the velocity at fixed normal is $fN+\mathrm{II}^{-1}\nabla_\partial f$. Differentiating
$\frac{d}{dt}\int_{K_t}U\dd x=\int_{\partial K_t}Uf\dd S$ and integrating
by parts on $\partial K$ yields the variation formula
\cite[Proposition~3.2]{KL21}
\[
 B(U)=\int_{\partial K}\bigl[(HU+\partial_NU)f^2
   -U\langle\mathrm{II}^{-1}\nabla_\partial f,\nabla_\partial f\rangle\bigr]\dd S.
\]
Thus $V'=0$ and $V''=B(1)$. For $p\ge2$, $k_p\in C^2(\R^n)$.
Write $U_p=k_p*\mathbf1_K$. Differentiating in both domains and using
symmetry in $x,y$, we obtain
\begin{align*}
 \frac{E_p'}2&=S_p:=\int_{\partial K}U_pf\dd S,\\
 \frac{E_p''}2&=B(U_p)+\int_{(\partial K)^2}k_p(x-y)f(x)f(y)\dd S_x\dd S_y .
\end{align*}

Write $A=D^2u$ and $w=\nabla_\partial u-\mathrm{II}^{-1}\nabla_\partial f$.
Since $\Delta u=0$,
\[
 \operatorname{div}\bigl((\nabla U\cdot v)v-UAv\bigr)
 =D^2U[v,v]-U|A|^2.
\]
On $\partial K$, $u_{NN}=-Hf-\Delta_\partial u$ and
$(AN)_{\rm tan}=\nabla_\partial f-\mathrm{II}\nabla_\partial u$.
The divergence theorem, integration by parts on $\partial K$ and completion
of the square give
\[
 B(U)=\int_K\bigl(D^2U[v,v]-U|A|^2\bigr)\dd x
                 -\int_{\partial K}U\mathrm{II}(w,w)\dd S .
\]
Thus $V''=-\mathcal H$, where
\[
 \mathcal H=\int_K|A|^2\dd x+\int_{\partial K}\mathrm{II}(w,w)\dd S.
\]
Since $\operatorname{div}v=0$ and $v\cdot N=f$, integration by parts yields
$S_p=\int_K\nabla U_p\cdot v\dd x$ and, in the two variables,
\[
 \int_{(\partial K)^2}k_p(x-y)f(x)f(y)\dd S_x\dd S_y
 =-\int_{K^2}D^2k_p(x-y)[v(x),v(y)]\dd x\dd y .
\]
Consequently,
\begin{align*}
 \frac{Q_p}2&=P_p-T_p-\frac{2(2n+p-1)}{2n+p}\frac{S_p^2}{E_p},\\
 P_p&=\frac12\int_{K^2}D^2k_p(x-y)[v(x)-v(y),v(x)-v(y)]\dd x\dd y,\\
 T_p&=\int_KU_p|A|^2\dd x+\int_{\partial K}U_p\mathrm{II}(w,w)\dd S .
\end{align*}

\vskip 3pt
\textbf{Step 2.} Estimate the second variation at $p=2$.

Let $c_{K_t}$ be the centroid of $K_t$. Write
\[
 M=\int_{\mathbb S^{n-1}}\theta\theta^\top\dd\nu(\theta),\qquad
 J(t)=\int_{K_t}(x-c_{K_t})^\top M(x-c_{K_t})\dd x.
\]
Then $\operatorname{tr}M=1$ and $E_2=2VJ$. For positive definite $M$,
a change of variables identifies $\sqrt{\det M}\,J(t)$ with the centered
polar moment of inertia of $M^{\frac12}K_t$. Thus $J^{\frac1{n+2}}$ is concave
by Hadwiger's inequality \cite[Section~2.3]{Hadwiger}. The singular case
follows by replacing $M$ with $M+\delta\,\mathrm{Id}$ and letting $\delta\downarrow0$.

Since $0<\frac{n+2}{2n+2}<1$, $J^{\frac1{2n+2}}$ is also concave.
Differentiating $E_2^{\frac1{2n+2}}=(2V)^{\frac1{2n+2}}J^{\frac1{2n+2}}$
twice at $t=0$ and using $V'(0)=0$, we obtain
\[
 \frac{Q_2}{E_2}
 =\frac{V''}{V}+(2n+2)\frac{(J^{\frac1{2n+2}})''}{J^{\frac1{2n+2}}}
 \le\frac{V''}{V}.
\]
Since $B^n\subseteq K$ and $\operatorname{tr}M=1$,
\[
 J(0)\ge\int_{B^n}(x-c_K)^\top M(x-c_K)\dd x
 =\frac{\omega_n}{n+2}+\omega_n c_K^\top M c_K
 \ge\frac{\omega_n}{n+2}.
\]
Thus, by $V''=-\mathcal H$,
\[
 Q_2\le 2J(0)V''=-2J(0)\mathcal H\le-\frac{2\omega_n}{n+2}\mathcal H.
\]

\vskip 3pt
\textbf{Step 3.} Compare $Q_p$ with $Q_2$.

For $2\le p\le\frac94$, we prove $Q_p-Q_2\le C_n(p-2)\mathcal H$.
Since $\int_Kv\dd x=0$ and $\operatorname{diam}K\le2n$, it follows from
the Payne--Weinberger inequality \cite{PW60} applied to each component of $v$ that
\[
 \int_K |v|^2\dd x
 \le \frac{(\operatorname{diam}K)^2}{\pi^2}\int_K |D^2u|^2\dd x
 \le C_n\mathcal H .
\]
For $|s|\le2n$, $(|s|^p)''\le p(p-1)(2n)^{p-2}\le 2+C_n(p-2)$.
Integrating against the probability measure $\nu$, we obtain
$D^2k_p-D^2k_2\le C_n(p-2)\,\mathrm{Id}$ on $2nB^n$. Hence
\[
 P_p-P_2\le C_n(p-2)\int_{K^2}|v(x)-v(y)|^2\dd x\dd y
 \le C_n(p-2)\mathcal H .
\]

The mean value theorem in $p$ gives $\big\||s|^p-s^2\big\|_{C^1([-2n,2n])}
\le C_n(p-2)$, so $\|U_p-U_2\|_{C^1(K)}\le C_n(p-2)$.
By rotation invariance the energy of $B^n$ is the positive continuous function
$p\mapsto\int_{(B^n)^2}|x_1-y_1|^p\dd x\dd y$. Since $B^n\subseteq K\subseteq nB^n$,
we have $c_n\le E_p\le C_n$.
Using $E_p=\int_K U_p\dd x$, $S_p=\int_K\nabla U_p\cdot v\dd x$ and
the $L^2$ bound for $v$, we find $|E_p-E_2|\le C_n(p-2)$,
$|S_p-S_2|\le C_n(p-2)\sqrt{\mathcal H}$ and $|S_p|\le C_n\sqrt{\mathcal H}$.
Thus
\[
 |T_p-T_2|\le C_n(p-2)\mathcal H,\qquad
 \left|\frac{S_p^2}{E_p}-\frac{S_2^2}{E_2}\right|\le C_n(p-2)\mathcal H .
\]
The negative square term in Step~1 has increasing coefficient
$\frac{2(2n+p-1)}{2n+p}$. This, together with the preceding estimates and the bound from Step~2, yields
\[
 Q_p\le Q_2+C_n(p-2)\mathcal H
     \le\left(-\frac{2\omega_n}{n+2}+C_n(p-2)\right)\mathcal H ,
\]
where $C_n>0$ depends only on $n$. With
\[
 \eta_n=\min\left\{\frac14,\frac{\omega_n}{(n+2)C_n}\right\},
\]
we have $Q_r\le-\frac{\omega_n}{n+2}\mathcal H$ for $2<r\le2+\eta_n$.

\vskip 3pt
\textbf{Step 4.} Pass to arbitrary convex bodies and characterize equality.

Let $F=E_{q,\nu}^{\frac1{n+q-1}}$ and $W=V^{\frac1n}$. Using $E_r\le C_n$,
$V\ge\omega_n$ and $V''=-\mathcal H$ in the estimate from Step~3, we have
\[
 (n+q-1)\frac{F''}{F}=\frac{Q_r}{E_r}\le\kappa_n\frac{V''}{V}
\]
for some $\kappa_n>0$. For the normalized velocity, $V'=0$, and hence
$\frac{V''}{V}=n\frac{W''}{W}$. The Hessians of $F$ and $W$ annihilate dilation
and translation directions, and at $p=r$ the affine normalization rescales each
functional. Thus for the original body and velocity,
\[
 (n+q-1)\frac{F''}{F}\le\kappa_n n\frac{W''}{W}.
\]

Along a smooth Minkowski segment, $W''\le0$, so $F$ is concave. In particular,
$F(K_t)\ge\min\{F(K),F(L)\}$. Since $K+L=K_t+K_{1-t}$ contains a translate of $K_t$,
we also have $W(K_t)\le W(K+L)$.
Thus $F''\le cW''$ with
\[
 c=\frac{\kappa_n n}{n+q-1}\frac{\min\{F(K),F(L)\}}{W(K+L)}>0.
\]
By midpoint concavity of $F-cW$, homogeneity and the Brunn--Minkowski
inequality for volume,
\[
 F(K+L)-F(K)-F(L)
 \ge c\bigl[W(K+L)-W(K)-W(L)\bigr]\ge0.
\]
By \eqref{eq:monotone-continuity}, $F$, $W$ and $c$ depend continuously on the
endpoints, so approximation by smooth bodies of positive curvature extends the
inequality to all $K,L\in\mathcal K^n$. Since $c>0$, equality for $F$ implies
equality for volume, hence homothety. The converse follows by homogeneity and
translation invariance.
\end{proof}

Theorem~\ref{aw:thm:angular-window} with $\dd\nu=\dd u$ and
\eqref{intro:eq:chord-energy} give the following inequality.

\begin{theorem}\label{thm:chord-bm-angular-window}
For $n+1<q\le n+3+\eta_n$ and $K,L\in\mathcal K^n$,
\[
 I_q(K+L)^{\frac1{n+q-1}}
 \ge I_q(K)^{\frac1{n+q-1}}+I_q(L)^{\frac1{n+q-1}},
\]
with equality if and only if $K$ and $L$ are homothetic.
\end{theorem}

\vskip 25pt
\section{\bfseries Applications}
\label{app:sec:chord-uniqueness}

\applicationsubsection{Uniqueness for chord measures}

For $q>0$, the chord measure \cite[Definition~4.1]{LXYZ} is defined on
Borel sets $\eta\subseteq\mathbb S^{n-1}$ by
\[
 F_q(K,\eta)=\frac{2q}{n\omega_n}
 \int_{\nu_K^{-1}(\eta)}
 \int_{S_x^+(K)}\rho_{K,x}(u)^{q-1}\dd u\,
                         \dd\mathcal H^{n-1}(x).
\]
Here $\rho_{K,x}(u)=\max\{a\ge0:x+au\in K\}$ and
$S_x^+(K)=\{u\in\mathbb S^{n-1}:\rho_{K,x}(u)>0\}$
\cite[Section~2]{LXYZ}. The map $\nu_K:\partial'K\to\mathbb S^{n-1}$ is the Gauss map of $K$,
defined on $\partial'K$, the set of points of $\partial K$ that have
a unique outer unit normal. Recall that
$\mathcal H^{n-1}(\partial K\setminus\partial'K)=0$
\cite[p.~84]{Schneider}.

For $q=1$, the chord measure is the surface area measure of $K$.
Lutwak, Xi, Yang, and Zhang~\cite[Theorem~1.1]{LXYZ} characterized the
measures that arise as $F_q(K,\cdot)$ for each $q>0$.

For $g\in C(\mathbb S^{n-1})$, write
\[
 [h_K+tg]=\{x\in\R^n:x\cdot u\le h_K(u)+tg(u)
                         \text{ for all }u\in \mathbb{S}^{n-1}\}.
\]
By \cite[Theorem~5.5]{LXYZ}, we have
\[
 \left.\frac{\dd}{\dd t}\right|_{t=0} I_q([h_K+tg])
       =\int_{\mathbb{S}^{n-1}}g(u)\dd F_q(K,u).
\]
For $g=h_K$, homogeneity gives
\begin{equation}\label{app:eq:chord-euler}
 \int_{\mathbb S^{n-1}}h_K(u)\dd F_q(K,u)=\alpha I_q(K),\qquad \alpha=n+q-1.
\end{equation}
For $0<q\le n+3+\eta_n$, Theorem~\ref{thm:chord-main}
implies that $G(t)=I_q((1-t)K+tL)^{\frac1\alpha}$ is concave on $[0,1]$.
By the variational formula with $g=h_L-h_K$ and \eqref{app:eq:chord-euler},
the bound $G'_+(0)\ge G(1)-G(0)$ becomes
\begin{equation}\label{app:eq:first-minkowski}
 \frac1\alpha\int_{\mathbb S^{n-1}}h_L(u)\dd F_q(K,u)
 \ge I_q(K)^{\frac{\alpha-1}{\alpha}}I_q(L)^{\frac1\alpha}.
\end{equation}
Equality holds precisely when $G$ is affine, that is, when $K$ and $L$ are homothetic.
Li, Wan, and Zhang~\cite[Theorem~8.1]{LWZ26} proved
\eqref{app:eq:first-minkowski} with its equality condition for
$q\in[1,n+2]\cup\{n+3\}$ in every dimension and for $q\in[1,5]$ in
$\mathbb{R}^2$.

\begin{theorem}\label{app:cor:chord-uniqueness}
Let $K,L\in\mathcal K^n$ and $0<q\le n+3+\eta_n$. If
$F_q(K,\cdot)=F_q(L,\cdot)$, then $K$ and $L$ are translates.
\end{theorem}
\begin{proof}
Since $F_q(K,\cdot)=F_q(L,\cdot)$, \eqref{app:eq:chord-euler} and
\eqref{app:eq:first-minkowski} give
\[
\begin{aligned}
 I_q(L)&=\frac{1}{\alpha}\int_{\mathbb S^{n-1}}h_L(u)\dd F_q(K,u)
       \ge I_q(K)^{\frac{\alpha-1}{\alpha}}I_q(L)^{\frac{1}{\alpha}},\\
 I_q(K)&=\frac{1}{\alpha}\int_{\mathbb S^{n-1}}h_K(u)\dd F_q(L,u)
       \ge I_q(L)^{\frac{\alpha-1}{\alpha}}I_q(K)^{\frac{1}{\alpha}}.
\end{aligned}
\]
Since $\alpha>1$, the two inequalities imply $I_q(K)=I_q(L)$,
hence equality in \eqref{app:eq:first-minkowski}.
Thus $L=x+\lambda K$ for some $x\in\R^n$ and $\lambda>0$. By homogeneity
and translation invariance, $\lambda^\alpha=\frac{I_q(L)}{I_q(K)}=1$,
and hence $L=x+K$.
\end{proof}

\applicationsubsection{A Urysohn inequality for fractional perimeters}

The mean width of $K\in\mathcal K^n$ is
\[
 w(K)=\frac1{n\omega_n}\int_{\mathbb S^{n-1}}
             \bigl(h_K(u)+h_K(-u)\bigr)\dd u.
\]

\begin{corollary}\label{app:cor:fractional-urysohn}
For $0<s<1$ and $K\in\mathcal K^n$,
\[
 P_s(K)\le P_s(B^n)\left(\frac{w(K)}2\right)^{n-s},
\]
with equality if and only if $K$ is a Euclidean ball.
\end{corollary}
\begin{proof}
For $q=1-s$, rotation invariance and \eqref{app:eq:chord-euler} give
$\dd F_q(B^n,u)=\frac{\alpha I_q(B^n)}{n\omega_n}\dd u$, where $\alpha=n-s$.
Since $I_{1-s}=\frac{s(1-s)}{n\omega_n}P_s$, \eqref{app:eq:first-minkowski} gives
\[
 \left(\frac{P_s(K)}{P_s(B^n)}\right)^{\frac1\alpha}
 =\left(\frac{I_q(K)}{I_q(B^n)}\right)^{\frac1\alpha}
 \le\frac1{n\omega_n}\int_{\mathbb S^{n-1}}h_K(u)\dd u
 =\frac{w(K)}2.
\]
Equality holds precisely when $K$ is homothetic to $B^n$.
\end{proof}

\vskip 25pt
\bibliographystyle{amsplain}

\begin{thebibliography}{99}
\vskip 5pt

\bibitem{BL76}
H. Brascamp, E. Lieb, On extensions of the Brunn--Minkowski and
Pr\'ekopa--Leindler theorems, including inequalities for log concave
functions, and with an application to the diffusion equation,
J. Funct. Anal. 22 (1976) 366-389.

\bibitem{Gardner02}
R. Gardner, The Brunn--Minkowski inequality,
Bull. Amer. Math. Soc. 39 (2002) 355-405.

\bibitem{Hadwiger50}
H. Hadwiger, Neue Integralrelationen f\"ur Eik\"orperpaare,
Acta Sci. Math. (Szeged) 13 (1950) 252--257.

\bibitem{Hadwiger}
H. Hadwiger, Konkave Eik\"orperfunktionale und h\"ohere
Tr\"agheitsmomente, Comment. Math. Helv. 30 (1956) 285-296.

\bibitem{HYZ25}
Y. Huang, D. Yang, G. Zhang, Minkowski problems for geometric measures,
Bull. Amer. Math. Soc. 62 (2025) 359-425.

\bibitem{KL21}
A. Kolesnikov, G. Livshyts, On the Gardner--Zvavitch conjecture:
symmetry in inequalities of Brunn--Minkowski type,
Adv. Math. 384 (2021) 107689.

\bibitem{LiTalk26}
Y. Li, Some recent progress on chord Minkowski problems,
lecture announcement, Tianyuan Mathematical Center in Southeast China,
16 October 2026, \url{https://tianyuan.xmu.edu.cn/info/1631/32721.htm}.

\bibitem{LWZ26}
H. Li, Y. Wan, Y. Zhang, Brunn--Minkowski inequalities for chord integrals,
preprint, \href{https://arxiv.org/abs/2611.34385}{arXiv:2611.34385}, 2026.

\bibitem{Ludwig13}
M. Ludwig, Anisotropic fractional perimeters,
J. Differential Geom. 96 (2014) 77-93.

\bibitem{LXYZ}
E. Lutwak, D. Xi, D. Yang, G. Zhang,
Chord measures in integral geometry and their Minkowski problems,
Comm. Pure Appl. Math. 77 (2024) 3277-3330.

\bibitem{McCann}
R. McCann, A convexity principle for interacting gases,
Adv. Math. 128 (1997) 153-179.

\bibitem{Moseeva22}
T. Moseeva, Integral identities for the boundary of a convex body,
Zap. Nauchn. Sem. POMI 510 (2022) 172-188 (in Russian).

\bibitem{PW60}
L. Payne, H. Weinberger, An optimal Poincar\'e inequality for convex domains,
Arch. Rational Mech. Anal. 5 (1960) 286-292.

\bibitem{Qin23}
L. Qin, Nonlocal energies of convex bodies and their log-Minkowski problems,
Adv. Math. 427 (2023) 109132.

\bibitem{SZ}
S. Sadovsky, G. Zhang, Brunn--Minkowski and reverse isoperimetric
inequalities for dual quermassintegrals,
Adv. Math. 480 (2025) 110456.

\bibitem{Santalo86}
L. A. Santal\'o, On the measure of line segments entirely contained in a convex body,
in: Aspects of Mathematics and its Applications,
North-Holland Math. Library 34 (1986) 677--687.

\bibitem{Santalo}
L. A. Santal\'o, Integral Geometry and Geometric Probability,
2nd ed., Cambridge University Press, Cambridge, 2004.

\bibitem{Schneider}
R. Schneider, Convex Bodies: The Brunn--Minkowski Theory,
second expanded edition, Cambridge University Press, Cambridge, 2014.

\bibitem{TXX23}
J. Tao, G. Xiong, J. Xiong, The logarithmic Minkowski inequality for cylinders,
Proc. Amer. Math. Soc. 151 (2023) 2143-2154.

\bibitem{XS08}
G. Xiong, X. Song, Inequalities for chord power integrals,
J. Korean Math. Soc. 45 (2008) 587-596.

\bibitem{XX22}
G. Xiong, J. Xiong, The Orlicz Minkowski problem for the electrostatic
$\mathfrak p$-capacity, Adv. Appl. Math. 137 (2022) 102339.

\bibitem{Zhang99}
G. Zhang, Dual kinematic formulas,
Trans. Amer. Math. Soc. 351 (1999) 985-995.

\bibitem{ZX18}
D. Zou, G. Xiong, A unified treatment for $L_p$ Brunn--Minkowski
type inequalities, Comm. Anal. Geom. 26 (2018) 435-460.

\end{thebibliography}

\end{document}